\documentstyle[12pt]{article}
\def\Bbb R{{\rm \bf R}}
\def\proclaim#1{\vskip2mm{\bf #1}\em}
\def\endproclaim{\em \vskip2mm}
\def\tag#1{\eqno(#1)}
\def\gathered{\begin{array}{c}}
\def\endgathered{\end{array}}
\def\text{\mbox}

\begin{document}

\title {On finding the surface admittance of an obstacle via the time domain enclosure method}
\author{Masaru IKEHATA\footnote{
Laboratory of Mathematics,
Graduate School of Engineering,
Hiroshima University,
Higashihiroshima 739-8527, JAPAN}}
\maketitle

\begin{abstract}
An inverse obstacle scattering problem for the electromagnetic wave
governed by the Maxwell system over a finite time interval is considered.
It is assumed that the wave satisfies the Leontovich boundary condition 
on the surface of an unknown obstacle.  The condition is described
by using an unknown positive function on the surface of the obstacle which is called the surface admittance.
The wave is generated at the initial time by a volumetric current source supported on a very small ball placed outside the obstacle
and only the electric component of the wave is observed on the same ball over a finite time interval.
It is shown that from the observed data one can extract information about the value of the surface admittance 
and the curvatures at the points on the surface nearest to the center of the ball.
This shows that a single shot contains a meaningful information about the quantitative state of the surface of the obstacle.

\noindent
AMS: 35R30, 35L50, 35Q61, 78A46, 78M35

\noindent KEY WORDS: enclosure method, inverse obstacle scattering problem, electromagnetic wave, obstacle,
Maxwell's equations, surface admittance, Leontovich boundary condition
\end{abstract}


\section{Introduction and statement of the results}

In this paper, we pursue further the possibility of the time domain enclosure method \cite{IR}
for the Maxwell system developed in \cite{IEM, IEM2}.
We consider an inverse obstacle scattering problem for the wave
governed by the Maxwell system in the time domain, in particular, over a {\it finite time interval}
unlike the time harmonic reduced case, see \cite{CK, KH, N}.

Let us formulate the problem more precisely.
We denote by $D$ the unknown obstacle.
We assume that $D$ is a non empty bounded open set of $\Bbb R^3$ with $C^2$-boundary such that $\Bbb R^3\setminus\overline D$
is connected.

We assume that the electric field $\mbox{\boldmath $E$}=\mbox{\boldmath $E$}(x,t)$ and magnetic field $\mbox{\boldmath $H$}
=\mbox{\boldmath $H$}(x,t)$
are generated only by the current density $\mbox{\boldmath $J$}=\mbox{\boldmath $J$}(x,t)$ at the initial time
located not far a way from the unknown obstacle. 
There should be several choices of current density
$\mbox{\boldmath $J$}$ as a model of antenna \cite{B, CB}. In this
paper, as considered in \cite{IEM, IEM2} we assume that $\mbox{\boldmath
$J$}$ takes the form
$$
\displaystyle
\mbox{\boldmath $J$}(x,t)=f(t)\chi_B(x)\mbox{\boldmath $a$},
\tag {1.1}
$$
where $\mbox{\boldmath $a$}$ is an arbitrary unit vector; $B$ is a (very small)
open ball satisfying $\overline B\cap\overline D=\emptyset$
and $\chi_B$ denotes the characteristic function of $B$; $f\in C^1[0,\,T]$ with $f(0)=0$.

Let $0<T<\infty$.
In this paper, we assume that

$\quad$

{\bf Assumption 1.}

(i)  the pair $(\mbox{\boldmath $E$}(t), \mbox{\boldmath $H$}(t))
\equiv (\mbox{\boldmath $E$}(\,\cdot\,,t), \mbox{\boldmath $H$}(\,\cdot\,,t))$ 
belongs to $C^1([0,\,T];L^2(\Bbb R^3\setminus\overline D)^3
\times L^2(\Bbb R^3\setminus\overline D)^3)$ as a function of $t$;

(ii)  for each $t\in[0,\,T]$, the pair $(\nabla\times\mbox{\boldmath $E$}(t),
\nabla\times\mbox{\boldmath $H$}(t))$ belongs to $L^2(\Bbb R^3\setminus\overline D)^3
\times L^2(\Bbb R^3\setminus\overline D)^3$;

(iii) it holds that
$$\displaystyle
\left\{
\begin{array}{l}
\displaystyle
\frac{d}{dt}\mbox{\boldmath $E$}
-\epsilon^{-1}\nabla\times\mbox{\boldmath $H$}
=
\epsilon^{-1}\mbox{\boldmath $J$},
\\
\\
\displaystyle
\frac{d}{dt}\mbox{\boldmath $H$}
+\mu^{-1}\nabla\times\mbox{\boldmath $E$}
=\mbox{\boldmath $0$},
\\
\\
\displaystyle
\mbox{\boldmath $E$}(0)=\mbox{\boldmath $0$},\\
\\
\displaystyle
\mbox{\boldmath $H$}(0)=\mbox{\boldmath $0$};
\end{array}
\right.
\tag {1.2}
$$

(iv) for each $t\in\,[0,\,T]$, $(\mbox{\boldmath $E$}(t), \mbox{\boldmath $H$}(t))$ satisfies,
in the sense of trace \cite{K} 
$$
\displaystyle
\begin{array}{ll}\mbox{\boldmath $\nu$}\times\mbox{\boldmath $H$}(t)
-\lambda\,
\mbox{\boldmath $\nu$}\times(\mbox{\boldmath $E$}(t)\times\mbox{\boldmath $\nu$})=\mbox{\boldmath $0$}
& \text{on}\,\partial D,
\end{array}
\tag {1.3}
$$
where $\lambda\in C^1(\partial D)$ and satisfies $\inf_{x\in\partial D}\lambda(x)>0$.

$\quad$

Note that $\mbox{\boldmath $\nu$}$ denotes the unit outward normal to $\partial D$.
The obstacle is embedded in
a medium like air (free space) which has constant electric permittivity $\epsilon(\,>0)$ and magnetic permeability $\mu(>0)$.
The boundary condition (1.3) is called the Leontovich boundary condition (\cite{AT, CK, KH, N}) and see
also \cite{KK} for the case when $\lambda$ is constant.
The quantity $1/\lambda$ is called the surface impedance, see \cite{AT} and thus $\lambda$ is called the admittance.
In what follows we use these equivalent forms without mentioning explicitly.
The existence of the admittance $\lambda$ causes the loss of the energy of the solution 
on the surface of the obstacle after stopping of the source supply.

In \cite{K} it is stated that the existence of $(\mbox{\boldmath $E$}(t),\mbox{\boldmath $H$}(t))$ satisfying (i)-(iv) can be
derived from the theory of $C_0$ contraction semigroups \cite{Y}.  However, therein the {\it detailed proof} is not given
as pointed out in \cite{IEM2}.
To make the logical relation clear, here we assume that our pair $(\mbox{\boldmath $E$}(t), \mbox{\boldmath $H$}(t))$
satisfies (i)-(iv).  This is our starting {\it assumption}.
It should be pointed out that Assumption ($\lambda$) in \cite{IEM2} implies the existence of such $(\mbox{\boldmath $E$}(t), \mbox{\boldmath $H$}(t))$
which ensures that conditions (i)-(iv) has a sense.

We consider the following problem.

$\quad$

{\bf\noindent Problem.}  Fix a large (to be determined later) $T<\infty$.  Observe $\mbox{\boldmath $E$}(t)$ on $B$
over the time interval $]0,\,T[$.
Extract information about the geometry of $D$ and the {\it values} of $\lambda$ on $\partial D$ from the observed data.

$\quad$

First of all let us recall the previous reslult on this problem.
Denote the solution of the system (1.2) in the case when $D=\emptyset$ by
$(\mbox{\boldmath $E$}_0(t), \mbox{\boldmath $H$}_0(t))$ with $\mbox{\boldmath $J$}$ given by (1.1).
Note that in this case, the solvabilty has been ensured by applying theory of $C^0$ contraction semigroups \cite{Y}.

Define the indicator function
$$\begin{array}{ll}
\displaystyle
I_{\mbox{\boldmath $J$}}(\tau,T)
=\int_B
\mbox{\boldmath $f$}(x,\tau)\cdot(\mbox{\boldmath $W$}_e-\mbox{\boldmath $V$}_e)dx,
&
\tau>0
\end{array}
\tag {1.4}
$$
where
$$
\begin{array}{ll}
\displaystyle
\mbox{\boldmath $W$}_e(x,\tau)
=\int_0^T e^{-\tau t}\mbox{\boldmath $E$}(x,t)dt,
&
\displaystyle
\mbox{\boldmath $V$}_e(x,\tau)
=\int_0^Te^{-\tau t}\mbox{\boldmath $E$}_0(x,t)dt
\end{array}
$$
and
$$\displaystyle
\mbox{\boldmath $f$}(x,\tau)=-\frac{\tau}{\epsilon}\int_0^Te^{-\tau t}\mbox{\boldmath $J$}(x,t)\,dt.
$$
And also, to describe another assumtion, we introduce here
$$\displaystyle
\mbox{\boldmath $W$}_m(x,\tau)
=\int_0^T e^{-\tau t}\mbox{\boldmath $H$}(x,t)dt.
$$
Using the same argument as that of \cite{IEM2} under Assumption 1, we know the following fatcts.

$\bullet$  The pair $(\mbox{\boldmath $W$}_e, \mbox{\boldmath $W$}_m)$ 
belongs to $L^2(\Bbb R^3\setminus\overline D)^3\times L^2(\Bbb R^3\setminus\overline D)^3$ with
$(\nabla\times\mbox{\boldmath $W$}_e, \nabla\times\mbox{\boldmath $W$}_m)
\in L^2(\Bbb R^3\setminus\overline D)^3\times L^2(\Bbb R^3\setminus\overline D)^3$.

$\bullet$  We have
$$
\left\{
\begin{array}{ll}
\displaystyle
\nabla\times\mbox{\boldmath $W$}_e+\tau\mu \mbox{\boldmath $W$}_m=-e^{-\tau T}\mu\mbox{\boldmath $H$}(x,T)
& \text{in}\,\Bbb R^3\setminus\overline D,\\
\\
\displaystyle
\nabla\times\mbox{\boldmath $W$}_m-\tau\epsilon \mbox{\boldmath $W$}_e
-\frac{\epsilon}{\tau}\mbox{\boldmath $f$}(x,\tau)
=e^{-\tau T}\epsilon \mbox{\boldmath $E$}(x,T)
& \text{in}\,\Bbb R^3\setminus\overline D.
\end{array}
\right.
\tag {1.5}
$$

$\bullet$  The boundary condition (1.3) remains valid in the sense of the trace \cite{K}
as mentioned above if $(\mbox{\boldmath $E$}(t),\mbox{\boldmath $H$}(t))$ 
is replaced with $(\mbox{\boldmath $W$}_e,\mbox{\boldmath $W$}_m)$.

$\bullet$  It holds that
$$\left\{
\begin{array}{ll}
\displaystyle
\nabla\cdot\mbox{\boldmath $W$}_m=0 & \text{in $\Bbb R^3\setminus\overline D$,}\\
\\
\displaystyle
\nabla\cdot\mbox{\boldmath $W$}_e=0 & \text{in $(\Bbb R^3\setminus\overline D)\setminus\overline B$.}\\
\end{array}
\right.
\tag {1.6}
$$

Note that, at this stage, each term on (1.3) does not have a point-wise meaning.  What we know is:
the left-hand side on (1.3) just belongs to the dual space of $H^{1/2}(\partial D)^3$.
In this paper, we introduce another assumption which states a {\it regularity} up to boundary.

$\quad$

{\bf\noindent Assumption 2.}
The functions $\mbox{\boldmath $W$}_e$ and $\mbox{\boldmath $W$}_m$ above belong
to $H^1$ on the intersection of an open neighbourhhod of $\partial D$ in $\Bbb R^3$ with $\Bbb R^3\setminus\overline D$.

$\quad$

This assumption makes us possible to treate vector-valued functions appeared in a dual paring {\it pointwise}.
Note that Assumption 2 is a special version of Assumption ($R$) introduced in \cite{IEM2} by virtue of (1.6).
However, for our purpose, it suffices to assume Assumption 2 instead of Assumption ($R$).
We believe that Assumption 2 should be removed.

Now, by Assumption 2, we have that both $\mbox{\boldmath $W$}_e$ and $\mbox{\boldmath $W$}_m$ 
belong to $H^1$ in the intersection of an open neighbourhood of $\partial D$ with $\Bbb R^3\setminus\overline D$.  Then, we see that
the boundary condition (1.3) is satisfied 
in the sense of the usual trace in $H^{1/2}(\partial D)^3$:
$$\begin{array}{ll}
\displaystyle
\mbox{\boldmath $\nu$}\times\mbox{\boldmath $W$}_m
-\lambda\,\mbox{\boldmath $\nu$}\times(\mbox{\boldmath $W$}_e\times\mbox{\boldmath $\nu$})=\mbox{\boldmath $0$}
& \text{on}\,\partial D.
\end{array}
$$
Note that this is equivalent to
$$\begin{array}{ll}
\displaystyle
\mbox{\boldmath $\nu$}\times(\mbox{\boldmath $W$}_m\times\mbox{\boldmath $\nu$})
+\lambda\,\mbox{\boldmath $\nu$}\times\mbox{\boldmath $W$}_e=\mbox{\boldmath $0$}
& \text{on}\,\partial D.
\end{array}
\tag {1.7}
$$
Moreover, note also that: since $\mbox{\boldmath $W$}_m\in L^2(\Bbb R^3\setminus\overline D)^3$ satifies
$\nabla\times\mbox{\boldmath $W$}_m\in L^2(\Bbb R^3\setminus\overline D)^3$,
from the first equation on (1.6) and by applying Corollary 1.1 on page 212 and Remark 2 on page 213 in \cite{DL3} one can conclude that 
$\mbox{\boldmath $W$}_m\in H^1(\Bbb R^3\setminus\overline D)^3$.

Set
$$
\displaystyle
\lambda_0=\sqrt{\frac{\epsilon}{\mu}}.
$$
As done in \cite{IEM}, we introduce two conditions (A.I) and (A.II) on $\lambda$ listed below:

$\quad$

(A.I) $\exists C>0$\,\,\, $\displaystyle\lambda(x)\ge\lambda_0+C$ for all $x\in\partial D$;

$\quad$

(A.II) $\exists C>0$\,$\exists C'>0$\,\,\, 
$\displaystyle C'\le\lambda(x)\le\lambda_0-C$ for all $x\in\partial D$.

$\quad$

Roughly speaking, we can say that: the condition (A.I)/(A.II) means that
the admittance $\lambda$ is greater/less than the special value 
$\lambda_0$ which is the admittance of free space \cite{AT}.

Define $\text{dist}\,(D,B)=\inf_{x\in D,\,y\in B}\,\vert x-y\vert$.

Under assumptions 1-2 we have already known that the following statement is true.

\noindent
\proclaim{\noindent Theorem 1.1(\cite{IEM2}).} 
Let $\mbox{\boldmath $a$}_j$,
$j=1,2$ be two linearly independent unit vectors. Let
$\mbox{\boldmath $J$}_j(x,t)=f(t)\chi_B(x)\mbox{\boldmath $a$}_j$
and $f$ satisfy
$$
\displaystyle
\exists\gamma\in\Bbb R\,\,\liminf_{\tau\rightarrow\infty}\tau^{\gamma}
\left\vert\int_0^Te^{-\tau t}f(t)\,dt\right\vert>0.
\tag {1.8}
$$

Then, we have:
$$
\displaystyle
\lim_{\tau\rightarrow\infty}e^{\tau T}\sum_{j=1}^2I_{\mbox{\boldmath $J$}_j}(\tau,T)
=
\left\{
\begin{array}{ll}
\displaystyle
0,
& \text{if $T\le2\sqrt{\mu\epsilon}\text{dist}\,(D,B)$,}\\
\\
\displaystyle
\infty,
&
\text{if $T>2\sqrt{\mu\epsilon}\text{dist}\,(D,B)$ and (A.I) is satified,}\\
\\
\displaystyle
-\infty,
&
\text{if $T>2\sqrt{\mu\epsilon}\text{dist}\,(D,B)$ and (A.II) is satisfied.}
\end{array}
\right.
$$

Moreover, if $\lambda$ satisfies (A.I) or (A.II), then
for all $T>2\sqrt{\mu\epsilon}\text{dist}\,(D,B)$
$$
\displaystyle
\lim_{\tau\rightarrow\infty}\frac{1}{\tau}
\log\left\vert
\sum_{j=1}^2I_{\mbox{\boldmath $J$}_j}(\tau,T)\right\vert
=-2\sqrt{\mu\epsilon}\text{dist}\,(D,B).
$$

\endproclaim

{\bf\noindent Remark 1.1.}
As described in Theorem 1.2 in \cite{IEM2}, all the statements of Theorem 1.1 are valid if $\mbox{\boldmath $V$}_e$ in $I_{\mbox{\boldmath $J$}}(\tau,T)$
is replaced with the unique weak solution $\mbox{\boldmath $V$}\in L^2(\Bbb R^3)^3$ with $\nabla\times
\mbox{\boldmath $V$}\in L^2(\Bbb R^3)^3$
of
$$\begin{array}{ll}
\displaystyle
\frac{1}{\mu\epsilon}\nabla\times\nabla\times\mbox{\boldmath $V$}
+\tau^2\mbox{\boldmath $V$}
+\mbox{\boldmath $f$}(x,\tau)=\mbox{\boldmath $0$}
& \text{in $\Bbb R^3$.}
\end{array}
\tag {1.9}
$$ 
In what follows, we denote by $\mbox{\boldmath $V$}_e^0$ the weak solution.
Roughly speaking, the reason why such a replacement is possible is the following.
Introduce another indicator function by the formula
$$
\displaystyle
\tilde{I}_{\mbox{\boldmath $f$}}(\tau,T)
=\int_B\mbox{\boldmath $f$}(x,\tau)\cdot(\mbox{\boldmath $W$}_e-\mbox{\boldmath $V$}_e^0)\,dx.
\tag {1.10}
$$
Using the simple facts
$$\displaystyle
\Vert\mbox{\boldmath $V$}_e-\mbox{\boldmath $V$}_e^0\Vert_{L^2(\Bbb R^3\setminus\overline D)}
=O(\tau^{-1}e^{-\tau T})
$$
and
$$
\displaystyle
\Vert\mbox{\boldmath $f$}\Vert_{L^2(B)}=O(\tau^{-1/2}),
$$
one has
$$\displaystyle
I_{\mbox{\boldmath $J$}}(\tau,T)
=\tilde{I}_{\mbox{\boldmath $J$}}(\tau,T)
+O(\tau^{-3/2}e^{-\tau T}).
\tag {1.11}
$$
Thus, one can transplant all the results in Theorem 1.1 into the case when the indicator function is given by (1.10).
This version's advantage is: no need of time domain computation of $\mbox{\boldmath $E$}_0$ in $\mbox{\boldmath $V$}_e$.

{\bf\noindent Remark 1.2.}  From Theorem 1.1. one can obtain another formula which has a similarity to
the original version of the enclosure method developed in \cite{I1}.  See (15) in \cite{IEM2}.

The main purpose of this paper is to go further beyond Theorem 1.1 under Assumptions 1 and 2.
Especially, we consider how to extract {\it quantitative} information about the state of the surafce of an unknown
obstacle using the time domain enclosure method.
For the purpose, we clarify the leading profile of the indicator functions (1.4) or (1.10) as $\tau\longrightarrow\infty$.

In what follows, we denote by $B_r(x)$ the open ball centered at $x$ with radius $r$.
Set $d_{\partial D}(p)=\inf_{y\in\partial D}\,\vert y-p\vert$ and
$\Lambda_{\partial D}(p)=\{y\in\partial D\,\vert\,\vert y-p\vert=d_{\partial D}(p)\}$.
To describe the formula, we recall some notion
in differential geometry.
Let $q\in\Lambda_{\partial D}(p)$.
Let $S_q(\partial D)$ and $S_q(\partial
B_{d_{\partial D}(p)}(p))$ denote the {\it shape operators} (or {\it Weingarten maps}) at $q$
of $\partial D$ and $\partial B_{d_{\partial D}(p)}(p)$ with
respect to $\nu_q$ and $-\nu_q$, respectively
(see \cite{O} for the notion of the shape operator).
Because $q$ attains the minimum of the function: $\partial D\ni y\longmapsto
\vert y-p\vert$, we have always $S_q(\partial B_{d_{\partial D}(p)}(p))-S_{q}(\partial D)\ge 0$
as the quadratic form on the common tangent space at $q$.

Now we are ready to state the main result in this paper.

\proclaim{\noindent Theorem 1.2.}
Assume that $\partial D$ is $C^4$ and $\lambda\in C^1(\partial D)$.
Assume that $\lambda$ satisfies (A1) or (A2).
Let $f$ satisfy (1.8) and $T>2\sqrt{\mu\epsilon}\,\text{dist}\,(D,B)$.
Assume that the set $\Lambda_{\partial D}(p)$ consists of finite points
and
$$\displaystyle
\text{det}\,(S_q(\partial B_{d_{\partial D}(p)}(p))-S_{q}(\partial D))>0\,\,\,
\forall q\in\Lambda_{\partial D}(p).
\tag {1.12}
$$
And also assume that $\mbox{\boldmath $\nu$}_q\times\mbox{\boldmath $a$}\not=\mbox{\boldmath $0$}$
for some $q\in\Lambda_{\partial D}(p)$.
Then, we have
$$\begin{array}{l}
\displaystyle
\,\,\,\,\,\,
\lim_{\tau\longrightarrow\infty}\tau^2e^{2\tau\,\sqrt{\mu\epsilon}\,\text{dist}\,(D,B)}\frac{\tilde{I}_{\mbox{\boldmath $J$}}(\tau,T)}
{\tilde{f}(\tau)^2}\\
\\
\displaystyle
=
\frac{\pi}{2}
\left(\frac{\eta}{d_{\partial D}(p)}\right)^2\frac{\lambda_0^2}{\epsilon^4}
\,
\sum_{q\in\Lambda_{\partial D}(p)}
k_q(p)\frac{\displaystyle
\lambda(q)-\lambda_0}
{\displaystyle
\lambda(q)+\lambda_0}
\vert\mbox{\boldmath $\nu$}_q\times\mbox{\boldmath $a$}\vert^2,
\end{array}
\tag {1.13}
$$
where
$$\displaystyle
\tilde{f}(\tau)=\int_0^T e^{-\tau t}f(t)\,dt
$$
and
$$\displaystyle
k_q(p)=\frac{1}{\displaystyle
\sqrt{\text{det}\,(S_q(\partial B_{d_{\partial D}(p)}(p))-S_{q}(\partial D))}}.
$$

\endproclaim

Note that $\vert\mbox{\boldmath $\nu$}_q\times(\mbox{\boldmath $a$}\times\mbox{\boldmath $\nu$}_q)\vert^2
=\vert\mbox{\boldmath $\nu$}_q\times\mbox{\boldmath $a$}\vert^2=1-(\mbox{\boldmath $\nu$}_q\cdot\mbox{\boldmath $a$})^2$.

Once we have the formula (1.13), as done in \cite{IC} for the scalar wave equation case,
we immediately obtain the following corollary.
To indicate the dependence of the indicator function on the surface admittance we write
$$\displaystyle
\tilde{I}_{\mbox{\boldmath $J$}}(\tau,T)=\tilde{I}_{\mbox{\boldmath $J$}}(\tau,T;\lambda).
$$

\proclaim{\noindent Corollary 1.1.}
Assume that $\partial D$ is $C^4$.
Let $\lambda_1$ and $\lambda_2$ belong to $C^1(\partial D)$ and satisfy (A1) or (A2).
Let $f$ satisfy (1.8) and $T>2\sqrt{\mu\epsilon}\,\text{dist}\,(D,B)$.
Assume that the set $\Lambda_{\partial D}(p)$ consists of finite points
and satisfies (1.12).
And also assume that $\mbox{\boldmath $\nu$}_q\times\mbox{\boldmath $a$}\not=\mbox{\boldmath $0$}$
for some $q\in\Lambda_{\partial D}(p)$.
Then, we have
$$\displaystyle
\lim_{\tau\longrightarrow\infty}
\frac{\tilde{I}_{\mbox{\boldmath $J$}}(\tau,T;\lambda_2)}
{\tilde{I}_{\mbox{\boldmath $J$}}(\tau,T;\lambda_1)}
=\frac{\displaystyle
\sum_{q\in\Lambda_{\partial D}(p)}
k_q(p)\frac{\displaystyle
\lambda_2(q)-\lambda_0}
{\displaystyle
\lambda_2(q)+\lambda_0}\vert\mbox{\boldmath $\nu$}_q\times\mbox{\boldmath $a$}\vert^2}
{\displaystyle
\sum_{q\in\Lambda_{\partial D}(p)}
k_q(p)\frac{\displaystyle
\lambda_1(q)-\lambda_0}
{\displaystyle
\lambda_1(q)+\lambda_0}\vert\mbox{\boldmath $\nu$}_q\times\mbox{\boldmath $a$}\vert^2}
$$
and its lower and upper estimates:
$$\displaystyle
\min_{q\in\Lambda_{\partial D}(p)}
\frac{
\displaystyle\frac{\lambda_2(q)-\lambda_0}{\lambda_2(q)+\lambda_0}
}
{
\displaystyle\frac{\lambda_1(q)-\lambda_0}{\lambda_1(q)+\lambda_0}
}
\le
\lim_{\tau\longrightarrow\infty}
\frac{\tilde{I}_{\mbox{\boldmath $J$}}(\tau,T;\lambda_2)}
{\tilde{I}_{\mbox{\boldmath $J$}}(\tau,T;\lambda_1)}
\le
\max_{q\in\Lambda_{\partial D}(p)}
\frac{
\displaystyle\frac{\lambda_2(q)-\lambda_0}{\lambda_2(q)+\lambda_0}
}
{
\displaystyle\frac{\lambda_1(q)-\lambda_0}{\lambda_1(q)+\lambda_0}
}.
\tag {1.14}
$$
In particulr, if $\Lambda_{\partial D}(p)$ consists of a single point $q\in\partial D$, we have
$$\displaystyle
\lim_{\tau\longrightarrow\infty}
\frac{\tilde{I}_{\mbox{\boldmath $J$}}(\tau,T;\lambda_2)}
{\tilde{I}_{\mbox{\boldmath $J$}}(\tau,T;\lambda_1)}
=\frac{
\displaystyle\frac{\lambda_2(q)-\lambda_0}{\lambda_2(q)+\lambda_0}
}
{
\displaystyle\frac{\lambda_1(q)-\lambda_0}{\lambda_1(q)+\lambda_0}
}.
\tag {1.15}
$$

\endproclaim

Estimates (1.14) and formula (1.15) are remarkable since they do not require information about the curvatures
of the surface of the obstacle in advance.
Note that if we know a point $q\in\Lambda_{\partial D}(p)$, then, all the intermediate points $p'$ on the segment connecting 
$q$ and $p$, satisfy $\Lambda_{\partial D}(p')=\{q\}$ and
$\displaystyle\text{det}\,(S_q(\partial B_{d_{\partial D}(p')}(p'))-S_q(\partial D))>0$.
Thus, one gets immediately the following corollary in which the set $\Lambda_{\partial D}(p)$ can be an {\it infinite} one, even, 
continuum.

\proclaim{\noindent Corollary 1.2.}
Assume that $\partial D$ is $C^4$.
Let $\lambda_1$ and $\lambda_2$ belong to $C^1(\partial D)$ and satisfy (A1) or (A2).
Le $p$ be an arbitrary point in $\Bbb R^3\setminus\overline D$ and $q\in\Lambda_{\partial D}(p)$.
Let $p'$ be an arbitrary point on the open segment $\{sq+(1-s)p\,\vert0<s<1\}$
and $B'$ an open ball centered at $p'$ satisfying $\overline B'\cap \overline D=\emptyset$.
Let $f$ satisfy (1.8) and $T>2\sqrt{\mu\epsilon}\,\text{dist}\,(D,B')$.
Let $\mbox{\boldmath $J$}'$ be the $\mbox{\boldmath $J$}$ given by (1.1) in which $B$ is replaced with $B'$.
And also assume that $\mbox{\boldmath $\nu$}_q\times\mbox{\boldmath $a$}\not=\mbox{\boldmath $0$}$.

Then, we have
$$\displaystyle
\lim_{\tau\longrightarrow\infty}
\frac{\tilde{I}_{\mbox{\boldmath $J$}'}(\tau,T;\lambda_2)}
{\tilde{I}_{\mbox{\boldmath $J$}'}(\tau,T;\lambda_1)}
=\frac{
\displaystyle\frac{\lambda_2(q)-\lambda_0}{\lambda_2(q)+\lambda_0}
}
{
\displaystyle\frac{\lambda_1(q)-\lambda_0}{\lambda_1(q)+\lambda_0}
}.
\tag {1.16}
$$

\endproclaim

Thus formula (1.16) can be used for {\it monitoring} of the quantitative state of the surface,
that is, the change of $\lambda_1$ to $\lambda_2$ of the surface admittance at a given monitoting point $q$ on the surface.

All the results mentioned above can be transplanted as follows.
\proclaim{\noindent Corollary 1.3.}
Theorem 1.2 and Corollaries 1.1-1.2 remain valid if $\tilde{I}_{\,\star\,}$ is replaced with
$I_{\,\star\,}$. 
\endproclaim

This can be seen as follows.  From (1.8) we have
$$\begin{array}{ll}
\displaystyle
\frac{\displaystyle
\tau^2e^{2\tau\sqrt{\mu\epsilon}\,\text{dist}\,(D,B)}}
{
\displaystyle
\tilde{f}(\tau)^2}
\tau^{-3/2}e^{-\tau T}
&
\displaystyle
=
\frac{\displaystyle
\tau^{2\gamma+1/2}e^{-\tau(T-2\sqrt{\mu\epsilon}\,\text{dist}\,(D,B))}}
{
\displaystyle
\tau^{2\gamma}\tilde{f}(\tau)^2}\\
\\
\displaystyle
&
\displaystyle
=O(\tau^{2\gamma+1/2}e^{-\tau(T-2\sqrt{\mu\epsilon}\,\text{dist}\,(D,B))}).
\end{array}
\tag {1.17}
$$
Thus if $T$ satisfies $T>2\sqrt{\mu\epsilon}\,\text{dist}\,(D,B)$, then
(1.11) and (1.17) ensure both quantities
$$\displaystyle
\frac{\displaystyle
\tau^2e^{2\tau\sqrt{\mu\epsilon}\,\text{dist}\,(D,B)}}
{
\displaystyle
\tilde{f}(\tau)^2}
\tilde{I}_{\,\star\,}
$$
and
$$\displaystyle
\frac{\displaystyle
\tau^2e^{2\tau\sqrt{\mu\epsilon}\,\text{dist}\,(D,B)}}
{
\displaystyle
\tilde{f}(\tau)^2}
I_{\,\star\,}
$$
have the same leading profile as $\tau\longrightarrow\infty$.

Finally, we show that Theorem 1.2 suggests us a procedure for
finding curvatures and $\lambda$ at an arbitrary point $q$ on $\Lambda_{\partial D}(p)$.
It is a translation of the procedure described in \cite{IC} in which the scalar wave equation
is considered.

$\quad$

{\bf\noindent Step 1.}  Choose three points $p_j$, $j=1,2,3$ on the segment connecting $p$ and $q$.
Denote by $B_j$ three open balls with very small radiuses centered at $p_j$ such that
$\overline B_1\cup\overline B_2\cup\overline B_3\subset\Bbb R^3\setminus\overline
D$. Note that we have $\Lambda_{\partial D}(p_j)=\{q\}$ and
$\displaystyle\text{det}\,(S_q(\partial B_{d_{\partial D}(p_j)}(p_j))-S_q(\partial D))>0$.

$\quad$

{\bf\noindent Step 2.}  Fix $T>2\,\max_{j}\,\sqrt{\mu\epsilon}\,\text{dist}\,(D, B_j)$
and generate  $\mbox{\boldmath $E$}$ and $\mbox{\boldmath $H$}$ on $B_j$ by the source $\mbox{\boldmath $J$}_j
=f(t)\chi_{B_j}\mbox{\boldmath $a$}$ for a fixed unit vector $\mbox{\boldmath $a$}$ with $\mbox{\boldmath $a$}\times\mbox{\boldmath $\nu$}_q
\not=\mbox{\boldmath $0$}$
and observe $\mbox{\boldmath $E$}$ on $B_j$ over the time interval $]0,\,T[$.

$\quad$

{\bf\noindent Step 3.}  Compute $\tilde{I}_{\mbox{\boldmath $J$}_j}(\tau,T)$ from the observation data in Step 2.

$\quad$

{\bf\noindent Step 4.}  Apply Theorem 1.2 to the case $B=B_j$.  Then, we have
$$\begin{array}{l}
\displaystyle
\lim_{\tau\longrightarrow\infty}
\tau^2e^{2\tau\,\sqrt{\mu\epsilon}\,\text{dist}\,(D,B_j)}
\frac{\displaystyle \tilde{I}_{\mbox{\boldmath $J$}_j}(\tau,T)}
{\tilde{f}(\tau)^2}
=\frac{\pi}{2}
\left(\frac{\eta}{d_{\partial D}(p_j)}\right)^2\frac{\lambda_0^4}{\epsilon^2}\vert\mbox{\boldmath $\nu$}_q\times\mbox{\boldmath $a$}\vert^2{\cal F}_j,
\end{array}
$$
where
$$\displaystyle
{\cal F}_j=\frac{\displaystyle\frac{\lambda(q)-\lambda_0}{\lambda(q)+\lambda_0}}
{\displaystyle\sqrt{\text{det}\,(S_q(\partial B_{d_{\partial D}(p_j)}(p_j))-S_q(\partial D))}},\,j=1,2,3.
$$

$\quad$

{\bf\noindent Step 5.}  Use the expression
$$\left\{
\begin{array}{l}
\displaystyle
\text{det}\,(S_q(\partial B_{d_{\partial D}(p_j)}(p_j))-S_q(\partial D))=
s_j^2-2H_{\partial D}(q)s_j+K_{\partial D}(q),\\
\\
\displaystyle
s_j=1/d_{\partial D}(p_j),
\end{array}
\right.
$$
where $H_{\partial D}(q)$ and $K_{\partial D}(q)$ denote the mean and Gauss curvatures at $q$ of $\partial D$
with respect to $\mbox{\boldmath $\nu$}_q$.
From ${\cal F}_j$ we have
$$
\displaystyle
\left(
\begin{array}{cc}
\displaystyle
\,\,\,
-(s_1{\cal F}_1^2-s_2{\cal F}_2^2)
&
\displaystyle
{\cal F}_1^2-{\cal F}_2^2\\
\\
\displaystyle
-(s_2{\cal F}_2^2-s_3{\cal F}_3^2)
&
\displaystyle
{\cal F}_2^2-{\cal F}_3^2
\end{array}
\right)
\left(\begin{array}{c}
\displaystyle
2 H_{\partial D}(q)
\\
\displaystyle
K_{\partial D}(q)
\end{array}
\right)
=
\left(\begin{array}{c}
\displaystyle
{\cal F}_2^2 s_2^2-{\cal F}_1^2 s_1^2
\\
\\
\displaystyle
{\cal F}_3^2 s_3^2-
{\cal F}_2^2 s_2^2
\end{array}
\right).
$$
Solving this linear system numerically, we may obtain $H_{\partial D}(q)$ and $K_{\partial D}(q)$.

$\quad$

{\bf\noindent Step 6.}  From ${\cal F}_j$ one has
$$\displaystyle
\left(\frac{\lambda(q)-\lambda_0}{\lambda(q)+\lambda_0}\right)^2=\frac{1}{3}
\sum_{j=1}^3{\cal F}_j^2(s_j^2-2H_{\partial D}(q)s_j+K_{\partial D}(q)).
$$

$\quad$

{\bf\noindent Step 7.}  From the signature of one of ${\cal F}_j$
one can know whether $\lambda(q)>\lambda_0$ or $\gamma(q)<\lambda_0$.

$\quad$

{\bf\noindent Step 8.}  From Steps 6 and 7 we may obtain $\lambda(q)$.

$\quad$

This paper is organized as follows.
In section 2, we give a proof of Theorem 1.2.
The proof is based on a rough asymptotic formula of the indicator function
as $\tau\longrightarrow\infty$ as stated in Lemma 2.1 which has been established in \cite{IEM2}.
The formula consists of two terms and remainder.
The treatement of the remainder is not a problem.
And the first term is explicitly given by (2.6)
as a Laplace type surface integral of $\mbox{\boldmath $V$}_e^0$ and its
rotation over $\partial D$.
Thus the key point is the profile of the second term
which is the energy integral of the so-called {\it reflected} solutions
given by (2.7).
Its asymptotic profile is stated as Theorem 2.1 which tells us
that the leading profle is also given as a Laplace type surface integral
of $\mbox{\boldmath $V$}_e^0$ and its rotation.
Then, using the leading profile of those two terms
which is described in Lemma 2.2 as an application of the Laplace method, 
we obtain the reading profile of the indicator function 
as stated in Theorem 1.2.

The proof of Theorem 2.1 is given in section 3.
First we construct a candidate of the {\it leading term} of the reflected solutions.
For the purpose we employ a combination of the {\it reflection principle} which has been established in \cite{IEM2}
and a cut-off argument in a neighbourhood of $\partial D$ with a cut-off parameter $\delta$.
Then the first and second terms of integral (2.7) is extracted as (3.7) in Lemma 3.1.  To show that the first term
is the reading profile we have to prove that the second term is small compared with first term.
We see that it is true if $\delta$ is properly chosen according to the size of $\tau$.
It's essence is described as Lemma 3.2.
The proof of Lemma 3.2 which is given in section 4, is a modification of 
the Lax-Phillips reflection argument \cite{LP} originally developed
for the study of the leading singularity of the {\it scattering kernel} 
for the scalar wave equation in the context of the Lax-Phillips scattering theory,
however, our version of the argument is rather straightforward.

\section{Proof of Theorem 1.2}

Define
$$\displaystyle
\mbox{\boldmath $V$}_m^0=-\frac{1}{\tau\mu}\,\nabla\times\mbox{\boldmath $V$}_e^0.
\tag {2.1}
$$
From this and (1.9) we have
$$
\left\{
\begin{array}{ll}
\displaystyle
\nabla\times\mbox{\boldmath $V$}_e^0+\tau\mu \mbox{\boldmath $V$}_m^0=\mbox{\boldmath $0$}
& \text{in}\,\Bbb R^3,\\
\\
\displaystyle
\nabla\times\mbox{\boldmath $V$}_m^0-\tau\epsilon \mbox{\boldmath $V$}_e^0
-\frac{\epsilon}{\tau}\mbox{\boldmath $f$}(x,\tau)
=\mbox{\boldmath $0$}
& \text{in}\,\Bbb R^3.
\end{array}
\right.
\tag {2.2}
$$
It is a due course to deduce that $\mbox{\boldmath $V$}_m^0\in H^1(\Bbb R^3)^3$ and $\mbox{\boldmath $V$}_e^0$ belongs to
$H^1$ in a neighbourhood of $\overline D$.

Define
$$
\left\{
\begin{array}{l}
\displaystyle
\mbox{\boldmath $R$}_e=\mbox{\boldmath $W$}_e-\mbox{\boldmath $V$}_e^0,
\\
\\
\displaystyle
\mbox{\boldmath $R$}_m=\mbox{\boldmath $W$}_m-\mbox{\boldmath $V$}_m^0.
\end{array}
\right.
$$
Note that $(\mbox{\boldmath $R$}_e,\mbox{\boldmath $R$}_m)$ belongs to $L^2(\Bbb R^3\setminus\overline D)^3\times
L^2(\Bbb R^3\setminus\overline D)^3$ with $(\nabla\times\mbox{\boldmath $R$}_e,\nabla\times\mbox{\boldmath $R$}_m)
\in\,L^2(\Bbb R^3\setminus\overline D)^3\times
L^2(\Bbb R^3\setminus\overline D)^3$;
$\mbox{\boldmath $R$}_m\in H^1(\Bbb R^3\setminus\overline D)^3$ and $\mbox{\boldmath $R$}_e$
belongs to $H^1$ in a neighbourhood of $\partial D$.

From (2.2) and (1.5) we see that $\mbox{\boldmath $R$}_e$ and $\mbox{\boldmath $R$}_m$ satisfy
$$
\left\{
\begin{array}{ll}
\displaystyle
\nabla\times\mbox{\boldmath $R$}_e+\tau\mu\mbox{\boldmath $R$}_m
=-e^{-\tau T}\mu\mbox{\boldmath $H$}(x,T)
& \text{in}\,\Bbb R^3\setminus\overline D,\\
\\
\displaystyle
\nabla\times\mbox{\boldmath $R$}_m-\tau\epsilon\mbox{\boldmath $R$}_e
=e^{-\tau T}\epsilon\mbox{\boldmath $E$}(x,T)
& \text{in}\,\Bbb R^3\setminus\overline D.
\end{array}
\right.
\tag {2.3}
$$
It follows from (1.7) that
$$\displaystyle
\mbox{\boldmath $\nu$}\times (\mbox{\boldmath $R$}_m\times\mbox{\boldmath $\nu$})
+\lambda\,\mbox{\boldmath $\nu$}\times\mbox{\boldmath $R$}_e
=-\mbox{\boldmath $\nu$}\times (\mbox{\boldmath $V$}_m^0\times\mbox{\boldmath $\nu$})
-\lambda\,\mbox{\boldmath $\nu$}\times\mbox{\boldmath $V$}_e^0.
\tag {2.4}
$$
From \cite{IEM2}, we have a rough asymptotic formula of the indicator function.

\proclaim{\noindent Lemma 2.1(\cite{IEM2}).}
We have, as $\tau\longrightarrow\infty$
$$
\displaystyle
\tilde{I}_{\mbox{\boldmath $f$}}(\tau,T)
=J(\tau)+E(\tau)+O(e^{-\tau T}\tau^{-1}),
\tag {2.5}
$$
where 
$$\displaystyle
J(\tau)
=\frac{1}{\mu\epsilon}\int_{\partial D}(\mbox{\boldmath $\nu$}\times\mbox{\boldmath $V$}_e^0)
\cdot\nabla\times\mbox{\boldmath $V$}_e^0\,dS-
\frac{\tau}{\epsilon}\int_{\partial D}\frac{1}{\lambda}\vert\mbox{\boldmath $V$}_m^0\times\mbox{\boldmath $\nu$}\vert^2dS
\tag {2.6}
$$
and
$$\displaystyle
E(\tau)
=\frac{\tau}{\epsilon}\left\{\int_{\Bbb R^3\setminus\overline D}(\tau\mu\vert\mbox{\boldmath $R$}_m\vert^2
+\tau\epsilon\vert\mbox{\boldmath $R$}_e\vert^2)\,dx
+\int_{\partial D}\frac{1}{\lambda}\vert\mbox{\boldmath $R$}_m\times\mbox{\boldmath $\nu$}\vert^2\,dS\right\}.
\tag {2.7}
$$

\endproclaim

Thus, the essential part of the proof of Theorem 1.2 should be the study of the asymptotic behaviour of
$J(\tau)$ and $E(\tau)$ as
$\tau\longrightarrow\infty$.  The asymptotic behaviour of $J(\tau)$
can be reduced to that of a Laplace-type integral \cite{BH}.
See \cite{IEM2}.
For that of $E(\tau)$, we have the following result, which enables us to make a reduction
of the study to a Laplace-type integral.

\proclaim{\noindent Theorem 2.1.}
Assume that $\partial D$ is $C^4$ an $\lambda\in C^2(\partial D)$.
Assume that $\lambda$ has a positive lower bound, the set $\Lambda_{\partial D}(p)$ consists of finite points,
and (1.12) is satisfied; there exists a point $q\in\Lambda_{\partial D}(p)$ such that $\lambda(q)\not=\lambda_0$
and that
$$\displaystyle
\mbox{\boldmath $\nu$}_q\times(\mbox{\boldmath $a$}\times\mbox{\boldmath $\nu$}_q)\not=\mbox{\boldmath $0$}.
\tag {2.8}
$$
Let $f$ satisfy (1.8) and
$$\displaystyle
T>\sqrt{\mu\epsilon}\,\text{dist}\,(D,B).
\tag {2.9}
$$
Then, we have
$$
\displaystyle
\lim_{\tau\longrightarrow\infty}\frac{E(\tau)}
{\displaystyle
J^*(\tau)
}=1,
\tag {2.10}
$$
where
$$\displaystyle
J^*(\tau)
=\frac{\tau}{\epsilon}
\int_{\partial D}
\frac{\displaystyle
\lambda-\lambda_0}
{\displaystyle
\lambda+\lambda_0}\,
(\mbox{\boldmath $\nu$}\times\mbox{\boldmath $V$}_m^0)\cdot
\mbox{\boldmath $V$}_{em}^0\,dS
\tag {2.11}
$$
and
$$\displaystyle
\mbox{\boldmath $V$}_{em}^0=\mbox{\boldmath $\nu$}\times(\mbox{\boldmath $V$}_e^0\times\mbox{\boldmath $\nu$})
-\frac{1}{\lambda}\,\mbox{\boldmath $\nu$}\times\mbox{\boldmath $V$}_m^0.
\tag {2.12}
$$
\endproclaim

The proof of Theorem 2.1 is given in Section 3.
Assumption (2.8) means that vector $\mbox{\boldmath $a$}$ is not parallel to $\mbox{\boldmath $\nu$}_q$ at $q$.
Note that the factor $2$ in the restriction $T>2\sqrt{\mu\epsilon}\,\text{dist}\,(D,B)$ in
Theorem 1.2 is dropped in (2.9).  The quantity $\sqrt{\mu\epsilon}\,\text{dist}\,(D,B)$
corresponds to the first arrival time of the wave generated at $t=0$ on $B$ and reached at $\partial D$ firstly.
The asymptotic formula (2.10) clarifies the effect on the leading profile of the energy of the reflected solutions
$\mbox{\boldmath $R$}_e$ and $\mbox{\boldmath $R$}_m$ in terms of the deviation of th surface admittance
from that of free-space admittance and the energy density of the incident wave.

To complete the proof of Theorem 1.2 we need the following asymptotic formulae
of $J(\tau)$ and $J^*(\tau)$ as $\tau\longrightarrow\infty$.

\proclaim{\noindent Lemma 2.2.}
We have
$$\begin{array}{l}
\displaystyle
\,\,\,\,\,\,
\lim_{\tau\longrightarrow\infty}
\tau^2 e^{2\tau\sqrt{\mu\epsilon}\,\text{dist}\,(D,B)}\frac{J(\tau)}{\displaystyle\tilde{f}(\tau)^2}\\
\\
\displaystyle
=\frac{\pi}{4}\left(\frac{\eta}{d_{\partial D}(p)}\right)^2\frac{\lambda_0^3}{\epsilon^4}
\sum_{q\in\Lambda_{\partial D}(p)}k_q(p)
\left(\frac{1}{\lambda_0}-\frac{1}{\lambda(q)}\right)
\vert\mbox{\boldmath $\nu$}_q\times(\mbox{\boldmath $a$}\times\mbox{\boldmath $\nu$}_q)\vert^2
\end{array}
\tag {2.13}
$$
and
$$\begin{array}{l}
\displaystyle
\,\,\,\,\,\,
\lim_{\tau\longrightarrow\infty}
\tau^2 e^{2\tau\sqrt{\mu\epsilon}\,\text{dist}\,(D,B)}\frac{J^*(\tau)}{\displaystyle\tilde{f}(\tau)^2}\\
\\
\displaystyle
=\frac{\pi}{4}\left(\frac{\eta}{d_{\partial D}(p)}\right)^2\frac{\lambda_0^3}{\epsilon^4}
\sum_{q\in\Lambda_{\partial D}(p)}k_q(p)
\frac{\displaystyle \lambda(q)-\lambda_0}
{\displaystyle \lambda(q)+\lambda_0}
\left(\frac{1}{\lambda_0}-\frac{1}{\lambda(q)}\right)
\vert\mbox{\boldmath $\nu$}_q\times(\mbox{\boldmath $a$}\times\mbox{\boldmath $\nu$}_q)\vert^2.
\end{array}
\tag {2.14}
$$

\endproclaim

{\it\noindent Proof.}
Using (2.1), (2.12) and a simple computation in vector analysis, one can rewrite the right-hand side on (2.6) as
$$
\,\,\,\,\,\,
\displaystyle
J(\tau)=\frac{\tau}{\epsilon}
\int_{\partial D}
\mbox{\boldmath $\nu$}\times\mbox{\boldmath $V$}_m^0\cdot\mbox{\boldmath $V$}_{em}^0\,dS.
\tag {2.15}
$$
Set
$$\displaystyle
v(x,\tau)=\frac{e^{-\tau\sqrt{\mu\epsilon}\,\vert x-p\vert}}{\vert x-p\vert}.
$$
By (18) in \cite{IEM} we have already shown that $\mbox{\boldmath $V$}_e^0$ has the form
$$\begin{array}{ll}
\displaystyle
\mbox{\boldmath $V$}_e^0
=K(\tau)\tilde{f}(\tau)v\mbox{\boldmath $M$}\mbox{\boldmath $a$}
& \text{in}\,\Bbb R^3\setminus\overline B,
\end{array}
\tag {2.16}
$$
where
$$
\left\{
\begin{array}{l}
\displaystyle
K(\tau)=\frac{\mu\tau\varphi(\tau\sqrt{\mu\epsilon}\,\eta)}{(\tau\sqrt{\mu\epsilon})^3},
\\
\\
\displaystyle
\varphi(\xi)=\xi\cosh\xi-\sinh\xi,
\end{array}
\right.
$$
$$
\left\{
\begin{array}{l}
\displaystyle
\mbox{\boldmath $M$}
=\mbox{\boldmath $M$}(x;\tau)
=AI_3-B\,\mbox{\boldmath $\omega$}_x\otimes
\mbox{\boldmath $\omega$}_x,
\\
\\
\displaystyle
A=A(x,\tau)=1+\frac{1}{\tau\sqrt{\mu\epsilon}}
\left(\frac{1}{\vert x-p\vert}
+\frac{1}{\tau\sqrt{\mu\epsilon}\vert x-p\vert^2}\right),\\
\\
\displaystyle
B=B(x,\tau)=1+\frac{3}{\tau\sqrt{\mu\epsilon}}
\left(\frac{1}{\vert x-p\vert}
+\frac{1}{\tau\sqrt{\mu\epsilon}\vert x-p\vert^2}\right)
\end{array}
\right.
$$
and
$$\displaystyle
\mbox{\boldmath $\omega$}_x=\frac{x-p}{\vert x-p\vert}.
$$
This yields
$$\displaystyle
\nabla\times\mbox{\boldmath $V$}_e^0
=-\tau\sqrt{\mu\epsilon}\,K(\tau)\tilde{f}(\tau)v
\,\left(1+\frac{1}{\tau\sqrt{\mu\epsilon}\vert x-p\vert}\right)
\mbox{\boldmath $\omega$}_x\times\mbox{\boldmath $a$}
$$
and thus (2.1) gives
$$\displaystyle
\mbox{\boldmath $V$}_m^0
=\lambda_0\,K(\tau)\tilde{f}(\tau)v
\,\left(1+\frac{1}{\tau\sqrt{\mu\epsilon}\vert x-p\vert}\right)
\mbox{\boldmath $\omega$}_x\times\mbox{\boldmath $a$}.
\tag {2.17}
$$
A combination of (2.16) and (2.17) gives
$$\begin{array}{l}
\displaystyle
\,\,\,\,\,\,
\mbox{\boldmath $V$}_{em}^0
\\
\\
\displaystyle
=
\mbox{\boldmath $\nu$}\times(\mbox{\boldmath $V$}_e^0\times\mbox{\boldmath $\nu$})
-\frac{1}{\lambda}
\,\mbox{\boldmath $\nu$}\times\mbox{\boldmath $V$}_m^0
\\
\\
\displaystyle
=\lambda_0K(\tau)\tilde{f}(\tau)v
\mbox{\boldmath $\nu$}\times
\left\{\frac{1}{\lambda_0}\,\mbox{\boldmath $M$}\mbox{\boldmath $a$}\times\mbox{\boldmath $\nu$}
-\frac{1}{\lambda}
\left(1+\frac{1}{\tau\sqrt{\mu\epsilon}\,\vert x-p\vert}\,\right)
\mbox{\boldmath $\omega$}_x\times\mbox{\boldmath $a$}\right\}
\\
\\
\displaystyle
=\lambda_0K(\tau)\tilde{f}(\tau)v
\mbox{\boldmath $\nu$}\times
\left({\cal D}(x)\mbox{\boldmath $a$}+O(\tau^{-1})\right),
\end{array}
\tag {2.18}
$$
where $O(\tau^{-1})$ means uniformly with respect to $x\in\partial D$ and
$$\displaystyle
{\cal D}(x)\mbox{\boldmath $a$}
=\frac{1}{\lambda_0}
\mbox{\boldmath $a$}\times\mbox{\boldmath $\nu$}
-\frac{1}{\lambda_0}
\,\left(\mbox{\boldmath $a$}\cdot\mbox{\boldmath $\omega$}_x\right)
\,\mbox{\boldmath $\omega$}_x\times\mbox{\boldmath $\nu$}
-\frac{1}{\lambda}\,
\mbox{\boldmath $\omega$}_x\times\mbox{\boldmath $a$}.
$$
Thus we obtain
$$
\begin{array}{l}
\,\,\,\,\,\,
\displaystyle
\mbox{\boldmath $\nu$}\times\mbox{\boldmath $V$}_m^0
\cdot\mbox{\boldmath $V$}_{em}^0\\
\\
\displaystyle
=\lambda_0^2K^2(\tau)\tilde{f}(\tau)^2 v^2
\left\{
\mbox{\boldmath $\nu$}\times\left(\mbox{\boldmath $\omega$}_x\times\mbox{\boldmath $a$}\right)
\cdot\mbox{\boldmath $\nu$}\times({\cal D}(x)\mbox{\boldmath $a$})
+O(\tau^{-1})\right\}
\end{array}
\tag {2.19}
$$
and (2.15) gives
$$\displaystyle
\frac{J(\tau)}{\displaystyle\tilde{f}(\tau)^2}
=\frac{\lambda_0^2}{\epsilon}\tau\,K^2(\tau)
\int_{\partial D}\left\{
\mbox{\boldmath $\nu$}\times\left(\mbox{\boldmath $\omega$}_x\times\mbox{\boldmath $a$}\right)
\cdot\mbox{\boldmath $\nu$}\times({\cal D}(x)\mbox{\boldmath $a$})
+O(\tau^{-1})\right\}v^2\,dS.
\tag {2.20}
$$
Note that if $x\in \Lambda_{\partial D}(p)$, then $\mbox{\boldmath $\nu$}$ at $x$ coincides with
$-\mbox{\boldmath $\omega$}_x$.  Thus we have
$$\displaystyle
{\cal D}(x)\mbox{\boldmath $a$}
=\left(\frac{1}{\lambda_0}-\frac{1}{\lambda}\right)\,\mbox{\boldmath $a$}\times\mbox{\boldmath $\nu$}
$$
and hence
$$\displaystyle
\mbox{\boldmath $\nu$}\times\left(\mbox{\boldmath $\omega$}_x\times\mbox{\boldmath $a$}\right)
\cdot\mbox{\boldmath $\nu$}\times({\cal D}(x)\mbox{\boldmath $a$})
=\left(\frac{1}{\lambda_0}-\frac{1}{\lambda}\right)
\vert\mbox{\boldmath $\nu$}\times(\mbox{\boldmath $a$}\times\mbox{\boldmath $\nu$})\vert^2.
\tag {2.21}
$$

It is well known that the Laplace method under the assumption
that $\Lambda_{\partial D}(p)$ is finite and satisfies (1.12), yields
$$
\displaystyle
\lim_{\tau\longrightarrow\infty}
\tau  e^{2\tau d_{\partial D}(p)}
\int_{\partial D}A(x)
\frac{e^{-2\tau\vert x-p\vert}}
{\vert x-p\vert^2}dS
=\frac{\pi}{d_{\partial D}(p)^2}
\sum_{q\in\Lambda_{\partial D}(p)}
k_q(p)A(q),
$$
where $A\in C^1(\partial D)$.  See \cite{BH}, for example.  The point is that the Hessian matrix of the function
$\partial D\ni x\longmapsto \vert x-p\vert$ at $q\in\Lambda_{\partial D}(p)$ is given by
the operator $S_q(\partial B_{d_{\partial D}(p)}(p))-S_q(\partial D)$.  See, for example, \cite{IEE} for this point.
Replacing $\tau$ above with $\tau\,\sqrt{\mu\epsilon}$, we obtain
$$\displaystyle
\lim_{\tau\longrightarrow\infty}
\tau  e^{2\tau\sqrt{\mu\epsilon}\, d_{\partial D}(p)}
\int_{\partial D}A(x)v^2dS
=\frac{\pi}{\sqrt{\mu\epsilon}\,d_{\partial D}(p)^2}
\sum_{q\in\Lambda_{\partial D}(p)}
k_q(p)A(q).
\tag {2.22}
$$
Note also that
$$\displaystyle
K(\tau)\sim
\tau^{-1}\frac{\displaystyle
\eta e^{\tau\eta\sqrt{\mu\epsilon}}}
{2\epsilon}
$$
and thus
$$\begin{array}{ll}
\displaystyle
\frac{\lambda_0^2}{\epsilon}\tau\,K^2(\tau)
&
\displaystyle
=\frac{\tau}{\mu}K^2(\tau)\\
\\
\displaystyle
&
\displaystyle
\sim
\frac{\tau}{\mu}
\left(\tau^{-1}\frac{\displaystyle
\eta e^{\tau\eta\sqrt{\mu\epsilon}}}
{2\epsilon}\right)^2\\
\\
\displaystyle
&
\displaystyle
=\frac{1}{\mu}
\left(\frac{\eta}{2\epsilon}\right)^2
\tau^{-1}e^{2\tau\eta\sqrt{\mu\epsilon}}.
\end{array}
\tag {2.23}
$$
Applying (2.22) to (2.20) and noting (2.21) and (2.23), we obtain
$$\begin{array}{l}
\,\,\,\,\,\,
\displaystyle
\frac{J(\tau)}{\displaystyle\tilde{f}(\tau)^2}\\
\\
\displaystyle
\sim
\frac{1}{\mu}
\left(\frac{\eta}{2\epsilon}\right)^2
\tau^{-1}e^{2\tau\eta\sqrt{\mu\epsilon}}
\int_{\partial D}\left\{
\mbox{\boldmath $\nu$}\times\left(\mbox{\boldmath $\omega$}_x\times\mbox{\boldmath $a$}\right)
\cdot\mbox{\boldmath $\nu$}\times({\cal D}(x)\mbox{\boldmath $a$})
+O(\tau^{-1})\right\}v^2\,dS
\\
\\
\displaystyle
\sim\frac{\pi}{4}\left(\frac{\eta}{d_{\partial D}(p)}\right)^2\frac{\lambda_0^3}{\epsilon^4}
e^{-2\tau\sqrt{\mu\epsilon}\,\text{dist}\,(D,B)}\tau^{-2}\sum_{q\in\Lambda_{\partial D}(p)}
k_q(p)\left(\frac{1}{\lambda_0}-\frac{1}{\lambda(q)}\right)
\vert\mbox{\boldmath $\nu$}_q\times(\mbox{\boldmath $a$}\times\mbox{\boldmath $\nu$}_q)\vert^2.
\end{array}
$$
This is nothing but (2.13).
Similary, from (2.11), (2.19), (2.21), (2.22) and (2.23) we obtain (2.14).

\noindent
$\Box$

Now applying (2.10), (2.13) and (2.14) to (2.5) and noting
$$\begin{array}{ll}
\displaystyle
\frac{\displaystyle
\tau^2e^{2\tau\sqrt{\mu\epsilon}\,\text{dist}\,(D,B)}}
{\tilde{f}(\tau)^2}\tau^{-1}e^{-\tau T}
&
\displaystyle
=
\frac{\displaystyle
\tau^{1+2\gamma}e^{-\tau(T-2\sqrt{\mu\epsilon}\,\text{dist}\,(D,B)\,)}}
{\tau^{2\gamma}\tilde{f}(\tau)^2}\\
\\
\displaystyle
&
\displaystyle
=O(\tau^{1+2\gamma} e^{-\tau(T-2\sqrt{\mu\epsilon}\,\text{dist}\,(D,B)\,)})
\end{array}
$$
provided (1.8) is satisfied,
we obtain (1.13).
This completes the proof of Theorem 1.1.

\section{Proof of Theorem 2.1}

We denote by $x^r$ the reflection across $\partial D$ of the point $x\in\Bbb R^3\setminus D$ with
$d_{\partial D}(x)<2\delta_0$ for a sufficiently small $\delta_0>0$.
It is given by $x^r=2q(x)-x$, where $q(x)$ denotes the unique point on $\partial D$ such that
$d_{\partial D}(x)=\vert x-q(x)\vert$.  Note that $q(x)$ is $C^2$ for $x\in\Bbb R^3\setminus D$ with $d_{\partial D}(x)<2\delta_0$ if $\partial D$ is $C^3$ (see \cite{GT}).
Define $\tilde{\lambda}(x)=\lambda(q(x))$ for $x\in\Bbb R^3\setminus D$ with $d_{\partial D}(x)<2\delta_0$.
The function $\tilde{\lambda}$ is $C^2$ therein and coincides with $\lambda(x)$ for $x\in\partial D$.

Choose a cutoff function $\phi_{\delta}\in C^{2}(\Bbb R^3)$ with
$0<\delta<\delta_0$ which satisfies
$0\le\phi_{\delta}(x)\le 1$; $\phi_{\delta}(x)=1$ if $d_{\partial D}(x)<\delta$;
$\phi_{\delta}(x)=0$ if $d_{\partial D}(x)>2\delta$;
$\vert\nabla\phi_{\delta}(x)\vert\le C\delta^{-1}$;
$\vert\nabla^2\phi_{\delta}(x)\vert\le C\delta^{-2}$.

Using the reflection across the boundary $\partial D$, in \cite{IEM} we have already constructed
from $\mbox{\boldmath $V$}_e^0$ in $D$ the vector field $(\mbox{\boldmath $V$}_e^0)^*$
for $x\in\Bbb R^3\setminus\overline D$ with $d_{\partial D}(x)<2\delta_0$
and another one 
$$
\displaystyle (\mbox{\boldmath $V$}_m^0)^*\equiv -\frac{1}{\tau\mu}\nabla\times\{(\mbox{\boldmath $V$}_e^0)^*\}
$$
which satisfy
$$\begin{array}{ll}
\displaystyle
(\mbox{\boldmath $V$}_e^0)^*=-\mbox{\boldmath $V$}_e^0 & \text{on $\partial D$}
\end{array}
$$
and
$$\begin{array}{ll}
\displaystyle
\mbox{\boldmath $\nu$}\times(\mbox{\boldmath $V$}_m^0)^*=
\mbox{\boldmath $\nu$}\times\mbox{\boldmath $V$}_m^0 & \text{on $\partial D$.}
\end{array}
\tag {3.1}
$$

Define
$$\displaystyle
\mbox{\boldmath $R$}_e^0
=\frac{\displaystyle\tilde{\lambda}-\lambda_0}
{\displaystyle\tilde{\lambda}+\lambda_0}
\phi_{\delta}
\,(\mbox{\boldmath $V$}_e^0)^*
$$
and
$$\displaystyle
\mbox{\boldmath $R$}_m^0
=\frac{\displaystyle\tilde{\lambda}-\lambda_0}
{\displaystyle\tilde{\lambda}+\lambda_0}
\,\phi_{\delta}\,(\mbox{\boldmath $V$}_m^0)^*.
\tag {3.2}
$$
The pair $(\mbox{\boldmath $R$}_e^0,\mbox{\boldmath $R$}_m^0)$ belongs to $H^1(\Bbb R^3\setminus\overline D)^3\times
H^1(\Bbb R^3\setminus\overline D)^3$ and depends on $\delta$.

Define
$$
\left\{
\begin{array}{l}
\displaystyle
\mbox{\boldmath $R$}_e^1=\mbox{\boldmath $R$}_e-\mbox{\boldmath $R$}_e^0,\\
\\
\displaystyle
\mbox{\boldmath $R$}_m^1=\mbox{\boldmath $R$}_m-\mbox{\boldmath $R$}_m^0.
\end{array}
\right.
$$
Since $\mbox{\boldmath $R$}_e$ and $\mbox{\boldmath $R$}_m$ satisfiy (2.4), 
we obtain
$$\begin{array}{l}
\,\,\,\,\,\,
\displaystyle
\mbox{\boldmath $\nu$}\times\mbox{\boldmath $R$}_m^1
-\lambda\,\mbox{\boldmath $\nu$}\times(\mbox{\boldmath $R$}_e^1\times\mbox{\boldmath $\nu$})\\
\\
\displaystyle
=\mbox{\boldmath $\nu$}\times\{\mbox{\boldmath $\nu$}\times(\mbox{\boldmath $R$}_m^1\times\mbox{\boldmath $\nu$})\}
-\lambda\,\mbox{\boldmath $\nu$}\times(\mbox{\boldmath $R$}_e^1\times\mbox{\boldmath $\nu$})\\
\\
\displaystyle
=-\lambda\,\mbox{\boldmath $\nu$}\times(\mbox{\boldmath $\nu$}\times\mbox{\boldmath $R$}_e^1)
-\mbox{\boldmath $\nu$}\times\left\{
\frac{2\lambda}{\lambda+\lambda_0}
\left\{\mbox{\boldmath $\nu$}\times(\mbox{\boldmath $V$}_m^0\times\mbox{\boldmath $\nu$})
+\lambda_0\,\mbox{\boldmath $\nu$}\times\mbox{\boldmath $V$}_e^0\right\}
\right\}
\\
\\
\displaystyle
\,\,\,
-\lambda\,\mbox{\boldmath $\nu$}\times(\mbox{\boldmath $R$}_e^1\times\mbox{\boldmath $\nu$})\\
\\
\displaystyle
=-\frac{2\lambda}{\lambda+\lambda_0}\left\{\mbox{\boldmath $\nu$}\times\mbox{\boldmath $V$}_m^0
-\lambda_0\,\mbox{\boldmath $\nu$}\times(\mbox{\boldmath $V$}_e^0\times\mbox{\boldmath $\nu$})\right\}.
\end{array}
\tag {3.3}
$$
Define
$$
\begin{array}{ll}
\displaystyle
\mbox{\boldmath $V$}_1
=\frac{2\lambda}{\lambda+\lambda_0}\left\{\mbox{\boldmath $\nu$}\times\mbox{\boldmath $V$}_m^0
-\lambda_0\,\mbox{\boldmath $\nu$}\times(\mbox{\boldmath $V$}_e^0\times\mbox{\boldmath $\nu$})\right\}
&
\text{on $\partial D$.}
\end{array}
\tag {3.4}
$$
It follows from (2.3) and (3.3) that $\mbox{\boldmath $R$}_e^1$ and $\mbox{\boldmath $R$}_m^1$ satisfy
$$
\left\{
\begin{array}{ll}
\displaystyle
\nabla\times\mbox{\boldmath $R$}_e^1+\tau\mu\mbox{\boldmath $R$}_m^1
=-(\nabla\times\mbox{\boldmath $R$}_e^0+\tau\mu\mbox{\boldmath $R$}_m^0)
-e^{-\tau T}\mu\mbox{\boldmath $H$}(x,T)
& \text{in}\,\Bbb R^3\setminus\overline D,\\
\\
\displaystyle
\nabla\times\mbox{\boldmath $R$}_m^1-\tau\epsilon\mbox{\boldmath $R$}_e^1
=-(\nabla\times\mbox{\boldmath $R$}_m^0-\tau\epsilon\mbox{\boldmath $R$}_e^0)
+e^{-\tau T}\epsilon\mbox{\boldmath $E$}(x,T)
& \text{in}\,\Bbb R^3\setminus\overline D
\end{array}
\right.
\tag {3.5}
$$
and
$$
\begin{array}{ll}
\displaystyle
\mbox{\boldmath $\nu$}\times\mbox{\boldmath $R$}_m^1
-\lambda\,\mbox{\boldmath $\nu$}\times(\mbox{\boldmath $R$}_e^1\times\mbox{\boldmath $\nu$})
=-\mbox{\boldmath $V$}_1 & \text{on}\,\partial D.
\end{array}
\tag {3.6}
$$

Now we are ready to state an asymptotic formula of $E(\tau)-J^*(\tau)$ as $\tau\longrightarrow\infty$
which extracts the main term involving $\mbox{\boldmath $\nu$}\times\mbox{\boldmath $R$}_m^1$
on $\partial D$.

\proclaim{\noindent Lemma 3.1.}
We have, as $\tau\longrightarrow\infty$
$$\begin{array}{ll}
\displaystyle
E(\tau)
&
\displaystyle
=J^*(\tau)
+\frac{\tau}{\epsilon}\int_{\partial D}
\mbox{\boldmath $V$}_{em}^0\cdot(\mbox{\boldmath $\nu$}\times\mbox{\boldmath $R$}_m^1)\,dS\\
\\
\displaystyle
&
\displaystyle
\,\,\,
+O(e^{-\tau T}(\tau^{-2}e^{-\tau\sqrt{\mu\epsilon}\,\text{dist}\,(D,B)}\vert\tilde{f}(\tau)\vert+\tau^{-1}e^{-\tau T})).
\end{array}
\tag {3.7}
$$

\endproclaim

{\it\noindent Proof.}
Recall (40) in \cite{IEM2}:
$$\begin{array}{l}
\displaystyle
\,\,\,\,\,\,
\int_{\Bbb R^3\setminus\overline D}(\tau\mu\vert\mbox{\boldmath $R$}_m\vert^2
+\tau\epsilon\vert\mbox{\boldmath $R$}_e\vert^2)\,dx
\\
\\
\displaystyle
\displaystyle
=\int_{\partial D}
\mbox{\boldmath $\nu$}\times(\mbox{\boldmath $R$}_m\times\mbox{\boldmath $\nu$})\cdot
(\mbox{\boldmath $\nu$}\times\mbox{\boldmath $R$}_e)\,dS
\\
\\
\displaystyle
\,\,\,
-e^{-\tau T}
\int_{\Bbb R^3\setminus\overline D}
(\mu\mbox{\boldmath $H$}(x,T)\cdot\mbox{\boldmath $R$}_m+
\epsilon\mbox{\boldmath $E$}(x,T)\cdot\mbox{\boldmath $R$}_e)\,dx.
\end{array}
$$
It follows from this and (2.7) that
$$\begin{array}{ll}
\displaystyle
E(\tau)
&
\displaystyle
=\frac{\tau}{\epsilon}
\int_{\partial D}
\left\{\mbox{\boldmath $\nu$}\times(\mbox{\boldmath $R$}_m\times\mbox{\boldmath $\nu$})\cdot
(\mbox{\boldmath $\nu$}\times\mbox{\boldmath $R$}_e)
+\frac{1}{\lambda}\vert\mbox{\boldmath $R$}_m\times\mbox{\boldmath $\nu$}\vert^2\right\}\,dS
\\
\\
\displaystyle
&
\,\,\,
\displaystyle
-e^{-\tau T}
\int_{\Bbb R^3\setminus\overline D}
(\mu\mbox{\boldmath $H$}(x,T)\cdot\mbox{\boldmath $R$}_m+
\epsilon\mbox{\boldmath $E$}(x,T)\cdot\mbox{\boldmath $R$}_e)\,dx.
\end{array}
\tag {3.8}
$$
From (2.4) we have
$$
\displaystyle
\mbox{\boldmath $\nu$}\times\mbox{\boldmath $R$}_e
=-\mbox{\boldmath $\nu$}\times\mbox{\boldmath $V$}_e^0
-\frac{1}{\lambda}\left\{\mbox{\boldmath $\nu$}\times(\mbox{\boldmath $V$}_m^0\times\mbox{\boldmath $\nu$})
+\mbox{\boldmath $\nu$}\times(\mbox{\boldmath $R$}_m\times\mbox{\boldmath $\nu$})\right\}.
$$
This gives
$$\begin{array}{l}
\displaystyle
\,\,\,\,\,\,
\displaystyle
-\mbox{\boldmath $\nu$}\times(\mbox{\boldmath $R$}_m\times\mbox{\boldmath $\nu$})\cdot(\mbox{\boldmath $\nu$}\times\mbox{\boldmath $R$}_e)\\
\\
\displaystyle
=\mbox{\boldmath $\nu$}\times(\mbox{\boldmath $R$}_m\times\mbox{\boldmath $\nu$})\cdot
\left\{\mbox{\boldmath $\nu$}\times\mbox{\boldmath $V$}_e^0
+
\frac{1}{\lambda}\left\{\mbox{\boldmath $\nu$}\times(\mbox{\boldmath $V$}_m^0\times\mbox{\boldmath $\nu$})
+\mbox{\boldmath $\nu$}\times(\mbox{\boldmath $R$}_m\times\mbox{\boldmath $\nu$})\right\}\right\}\\
\\
\displaystyle
=\frac{1}{\lambda}\vert\mbox{\boldmath $R$}_m\times\mbox{\boldmath $\nu$}\vert^2
+\mbox{\boldmath $V$}_{em}^0\cdot(\mbox{\boldmath $R$}_m\times\mbox{\boldmath $\nu$}).
\end{array}
$$
Substituting $\mbox{\boldmath $R$}_m=\mbox{\boldmath $R$}_m^0+\mbox{\boldmath $R$}_m^1$ into 
the second term on this right-hand side and using (3.1) and (3.2), we obtain 
$$\begin{array}{l}
\displaystyle
\,\,\,\,\,\,
\mbox{\boldmath $\nu$}\times(\mbox{\boldmath $R$}_m\times\mbox{\boldmath $\nu$})\cdot(\mbox{\boldmath $\nu$}\times\mbox{\boldmath $R$}_e)
+\frac{1}{\lambda}\vert\mbox{\boldmath $R$}_m\times\mbox{\boldmath $\nu$}\vert^2\\
\\
\displaystyle
=-\mbox{\boldmath $V$}_{em}^0\cdot(\mbox{\boldmath $R$}_m\times\mbox{\boldmath $\nu$})\\
\\
\displaystyle
=\frac{\displaystyle
\lambda-\lambda_0}
{\displaystyle
\lambda+\lambda_0}\,
\mbox{\boldmath $V$}_{em}^0\cdot(\mbox{\boldmath $\nu$}\times\mbox{\boldmath $V$}_m^0)
+\mbox{\boldmath $V$}_{em}^0\cdot(\mbox{\boldmath $\nu$}\times\mbox{\boldmath $R$}_m^1).
\end{array}
$$
Thus (3.8) becomes
$$\begin{array}{ll}
\displaystyle
E(\tau)
&
\displaystyle
=\frac{\tau}{\epsilon}\int_{\partial D}
\frac{\displaystyle
\lambda-\lambda_0}
{\displaystyle
\lambda+\lambda_0}\,
\mbox{\boldmath $V$}_{em}^0\cdot(\mbox{\boldmath $\nu$}\times\mbox{\boldmath $V$}_m^0)\,dS
+\frac{\tau}{\epsilon}\int_{\partial D}
\mbox{\boldmath $V$}_{em}^0\cdot(\mbox{\boldmath $\nu$}\times\mbox{\boldmath $R$}_m^1)\,dS\\
\\
\displaystyle
&
\displaystyle
\,\,\,
-e^{-\tau T}
\int_{\Bbb R^3\setminus\overline D}
(\mu\mbox{\boldmath $H$}(x,T)\cdot\mbox{\boldmath $R$}_m+
\epsilon\mbox{\boldmath $E$}(x,T)\cdot\mbox{\boldmath $R$}_e)\,dx.
\end{array}
\tag {3.9}
$$
By Lemma 3.2 in \cite{IEM2} we have
$$\displaystyle
\Vert\mbox{\boldmath $R$}_e\Vert_{L^2(\Bbb R^3\setminus\overline D)}
=\Vert\mbox{\boldmath $R$}_m\Vert_{L^2(\Bbb R^3\setminus\overline D)}
=O(\tau^{-1/2}\Vert\mbox{\boldmath $V$}_{em}^0\Vert_{L^2(\partial D)}
+\tau^{-1}e^{-\tau T}).
\tag {3.10}
$$
Here we make use of the following asymptotic formula which can be shown similarily as formulae in Lemma 2.2
by using (2.18):
$$\begin{array}{l}
\displaystyle
\,\,\,\,\,\,
\lim_{\tau\longrightarrow\infty}
\tau^3 e^{2\tau\sqrt{\mu\epsilon}\,\text{dist}\,(D,B)}\frac{\displaystyle
\int_{\partial D}\vert\mbox{\boldmath $V$}_{em}^0\vert^2\,dS
}{\displaystyle\tilde{f}(\tau)^2}\\
\\
\displaystyle
=\frac{\pi}{4}\left(\frac{\eta}{d_{\partial D}(p)}\right)^2\frac{\lambda_0^3}{\epsilon^3}
\sum_{q\in\Lambda_{\partial D}(p)}k_q(p)
\left(\frac{1}{\lambda_0}-\frac{1}{\lambda(q)}\right)^2
\vert\mbox{\boldmath $\nu$}_q\times(\mbox{\boldmath $a$}\times\mbox{\boldmath $\nu$}_q)\vert^2.
\end{array}
\tag {3.11}
$$
Applying this to the right-hand side on (3.10), we obtain
$$\displaystyle
\Vert\mbox{\boldmath $R$}_e\Vert_{L^2(\Bbb R^3\setminus\overline D)}
=\Vert\mbox{\boldmath $R$}_m\Vert_{L^2(\Bbb R^3\setminus\overline D)}
=O(\tau^{-2}e^{-\tau\sqrt{\mu\epsilon}\,\text{dist}\,(D,B)}\vert\tilde{f}(\tau)\vert
+\tau^{-1}e^{-\tau T}).
\tag {3.12}
$$
Now a combination of (3.9) and (3.12) yields (3.7).

\noindent
$\Box$

Thus, the problem is: clarify the asymptotic behaviour of
$\mbox{\boldmath $\nu$}\times\mbox{\boldmath $R$}_m^1$ on $\partial D$
as $\tau\longrightarrow\infty$.  The point is the choice of $\delta$.

\proclaim{\noindent Lemma 3.2.}
Choose $\delta=\tau^{-1/2}$.
We have
$$\displaystyle
\lim_{\tau\longrightarrow\infty}\tau^3\,e^{2\tau\,\sqrt{\mu\epsilon}\,\text{dist}\,(D,B)}
\frac{\displaystyle
\Vert\mbox{\boldmath $\nu$}\times\mbox{\boldmath $R$}_m^1\Vert_{L^2(\partial D)}^2}
{\tilde{f}(\tau)^2}=0.
\tag {3.13}
$$

\endproclaim

\noindent
The proof of Leema 3.2 is given in Section 4.

Now choose $\delta$ in the pair $(\mbox{\boldmath $R$}_e^0, \mbox{\boldmath $R$}_m^0)$ as
that of Lemma 3.2.

Write
$$\begin{array}{l}
\displaystyle
\,\,\,\,\,\,
\left\vert
\frac{\displaystyle
\frac{\tau}{\epsilon}
\int_{\partial D}\mbox{\boldmath $V$}_{em}^0\cdot(\mbox{\boldmath $\nu$}\times\mbox{\boldmath $R$}_m^1)\,dS}
{\displaystyle J^*(\tau)}
\right\vert
\\
\\
\displaystyle
\le
\frac{\displaystyle
\tau^{3/2}e^{\tau\,\sqrt{\mu\epsilon}\,\text{dist}\,(D,B)}
\left\Vert\mbox{\boldmath $V$}_{em}^0\right\Vert_{\partial D}
\cdot
\tau^{3/2}e^{\tau\,\sqrt{\mu\epsilon}\,\text{dist}\,(D,B)}
\Vert\mbox{\boldmath $\nu$}\times\mbox{\boldmath $R$}_m^1\Vert_{L^2(\partial D)}}
{\displaystyle
\epsilon
\tau^2
e^{2\tau\,\sqrt{\mu\epsilon}\,\text{dist}\,(D,B)}
\left\vert
J^*(\tau)
\right\vert
}.
\\
\\
\displaystyle
\end{array}
$$
Applying (2.14), (3.11) and (3.13) to this right-hand side, we obtain
$$
\displaystyle
\lim_{\tau\longrightarrow\infty}
\frac{\displaystyle
\frac{\tau}{\epsilon}
\int_{\partial D}\mbox{\boldmath $V$}_{em}^0\cdot(\mbox{\boldmath $\nu$}\times\mbox{\boldmath $R$}_m^1)\,dS}
{\displaystyle J^*(\tau)}
=0.
\tag {3.14}
$$
Write
$$\begin{array}{l}
\,\,\,\,\,\,
\displaystyle
\frac{O(e^{-\tau T}(\tau^{-2}e^{-\tau\sqrt{\mu\epsilon}\,\text{dist}\,(D,B)}\vert\tilde{f}(\tau)\vert+\tau^{-1}e^{-\tau T}))}
{J^*(\tau)}
\\
\\
\displaystyle
=\frac{
O(
e^{-\tau T}(e^{\tau\sqrt{\mu\epsilon}\,\text{dist}\,(D,B)}\vert\tilde{f}(\tau)\vert^{-1}+\tau e^{-\tau T}
e^{2\tau\sqrt{\mu\epsilon}\text{dist}\,(D,B)}
\vert\tilde{f}(\tau)\vert^{-2})
)
}
{\displaystyle
\tau^2e^{2\sqrt{\mu\epsilon}\,\text{dist}\,(D,B)}J^*(\tau)\tilde{f}(\tau)^{-2}}.
\end{array}
\tag {3.15}
$$
Note that, if $f$ and $T$ satisfiy (1.8) and (2.9), respectively,
then we have, as $\tau\longrightarrow\infty$
$$\begin{array}{l}
\displaystyle
\,\,\,\,\,\,
e^{-\tau T}(e^{\tau\sqrt{\mu\epsilon}\,\text{dist}\,(D,B)}\vert\tilde{f}(\tau)\vert^{-1}+\tau e^{-\tau T}
e^{2\tau\sqrt{\mu\epsilon}\text{dist}\,(D,B)}
\vert\tilde{f}(\tau)\vert^{-2})\\
\\
\displaystyle
=\tau^{\gamma}e^{-\tau\,(T-\sqrt{\mu\epsilon}\,\text{dist}\,(D,B))}\tau^{-\gamma}\vert\tilde{f}(\tau)\vert^{-1}
+\tau^{-1+2\gamma} e^{-2\tau\,(T-\sqrt{\mu\epsilon}\,\text{dist}\,(D,B))}\tau^{-2\gamma}\vert \tilde{f}(\tau)\vert^{-2}
\longrightarrow 0.
\end{array}
$$
Now, applying this to (3.15) with the help of (2.14), we see that the left-hand side on (3.15) converges to $0$
as $\tau\longrightarrow\infty$.  Applying this and (3.14) to the right-hand side on (3.7), we obtain (2.10).

\section{Proof of Lemma 3.2}

In this section, we denote by $C$ several positive constants independen of $\delta$ and $\tau$.

\proclaim{\noindent Lemma 4.1.}
We have
$$\begin{array}{l}
\displaystyle
\,\,\,\,\,\,
\Vert\mbox{\boldmath $\nu$}\times\mbox{\boldmath $R$}_m^1\Vert_{L^2(\partial D)}^2
\\
\\
\displaystyle
\le C(\Vert\mbox{\boldmath $V$}_1\Vert_{L^2(\partial D)}^2
+\tau^{-1}\Vert\mbox{\boldmath $F$}_1\Vert_{L^2(U_{\delta})}^2
+\tau^{-1}\Vert\mbox{\boldmath $F$}_2\Vert_{L^2(U_{\delta})}^2)+O(\tau^{-1}e^{-2\tau T}),
\end{array}
\tag {4.1}
$$
where $U_{\delta}=\{x\in\Bbb R^3\setminus\overline D\,\vert\,d_{\partial D}(x)<2\delta\}$
and
$$\left\{
\begin{array}{l}
\displaystyle
\mbox{\boldmath $F$}_1=\nabla\times\mbox{\boldmath $R$}_e^0+\tau\mu\mbox{\boldmath $R$}_m^0,\\
\\
\displaystyle
\mbox{\boldmath $F$}_2=\nabla\times\mbox{\boldmath $R$}_m^0-\tau\epsilon\mbox{\boldmath $R$}_e^0.
\end{array}
\right.
\tag {4.2}
$$

\endproclaim

{\it\noindent Proof.}
Taking the inner product of the both sides of the first equation on (3.5) with $\mbox{\boldmath $R$}_m^1$, we obtain
$$\displaystyle
(\nabla\times\mbox{\boldmath $R$}_e^1)\cdot\mbox{\boldmath $R$}_m^1+\tau\mu\vert\mbox{\boldmath $R$}_m^1\vert^2
=-\mbox{\boldmath $F$}_1
\cdot\mbox{\boldmath $R$}_m^1
-e^{-\tau T}\mu\mbox{\boldmath $H$}(x,T)\cdot\mbox{\boldmath $R$}_m^1.
\tag {4.3}
$$
Taking the inner product of the both sides of the second equation on (3.5) with $\mbox{\boldmath $R$}_e^1$, we obtain
$$
\displaystyle
(\nabla\times\mbox{\boldmath $R$}_m^1)\cdot\mbox{\boldmath $R$}_e^1-\tau\epsilon\vert\mbox{\boldmath $R$}_e^1\vert^2
=-
\mbox{\boldmath $F$}_2\cdot\mbox{\boldmath $R$}_e^1
+
e^{-\tau T}\epsilon\mbox{\boldmath $E$}(x,T)\cdot\mbox{\boldmath $R$}_e^1.
\tag {4.4}
$$
By virtue of the fact that $\mbox{\boldmath $R$}_m^1\in H^1(\Bbb R^3\setminus\overline D)^3$
and $\mbox{\boldmath $R$}_e^1\in H(\text{curl},\Bbb R^3\setminus\overline D)$, we have
$$
\displaystyle
\,\,\,\,\,\,
\int_{\Bbb R^3\setminus\overline D}\nabla\cdot(\mbox{\boldmath $R$}_e^1\times\mbox{\boldmath $R$}_m^1)\,dx
=<\mbox{\boldmath $R$}_e^1\times\mbox{\boldmath $\nu$}, \mbox{\boldmath $\nu$}\times(\mbox{\boldmath $R$}_m^1\times\mbox{\boldmath $\nu$})>_{-1/2,1/2},
$$
where this right-hand side denotes the value
of the bounded linear functional $\mbox{\boldmath $R$}_e^1\times\mbox{\boldmath $\nu$}$ on $H^{1/2}(\partial D)^3$
of $\mbox{\boldmath $\nu$}\times(\mbox{\boldmath $R$}_m^1\times\mbox{\boldmath $\nu$})\in H^{1/2}(\partial D)^3$.
However, $\mbox{\boldmath $R$}_e^1$ belongs to $H^1$ in a neighbourhood of $\partial D$ this coincides with the integral
$$
\displaystyle
-\int_{\partial D}
\mbox{\boldmath $\nu$}\times(\mbox{\boldmath $R$}_m^1\times\mbox{\boldmath $\nu$})\cdot
(\mbox{\boldmath $\nu$}\times\mbox{\boldmath $R$}_e^1)\,dS.
$$
Note also that
$$\displaystyle
\nabla\cdot(\mbox{\boldmath $R$}_e^1\times\mbox{\boldmath $R$}_m^1)
=(\nabla\times\mbox{\boldmath $R$}_e^1)\cdot\mbox{\boldmath $R$}_m^1
-(\nabla\times\mbox{\boldmath $R$}_m^1)\cdot\mbox{\boldmath $R$}_e^1.
$$
From these, (4.2), (4.3) and (4.4) we obtain
$$\begin{array}{l}
\displaystyle
\,\,\,\,\,\,
\int_{\Bbb R^3\setminus\overline D}(\tau\mu\vert\mbox{\boldmath $R$}_m^1\vert^2
+\tau\epsilon\vert\mbox{\boldmath $R$}_e^1\vert^2)\,dx
\\
\\
\displaystyle
\displaystyle
=\int_{\partial D}
\mbox{\boldmath $\nu$}\times(\mbox{\boldmath $R$}_m^1\times\mbox{\boldmath $\nu$})\cdot
(\mbox{\boldmath $\nu$}\times\mbox{\boldmath $R$}_e^1)\,dS
\\
\\
\displaystyle
\,\,\,
-
\int_{\Bbb R^3\setminus\overline D}
\mbox{\boldmath $F$}_1\cdot\mbox{\boldmath $R$}_m^1\,dx
+\int_{\Bbb R^3\setminus\overline D}
\mbox{\boldmath $F$}_2\cdot\mbox{\boldmath $R$}_e^1\,dx
\\
\\
\displaystyle
\,\,\,
-e^{-\tau T}
\int_{\Bbb R^3\setminus\overline D}
(\mu\mbox{\boldmath $H$}(x,T)\cdot\mbox{\boldmath $R$}_m^1+
\epsilon\mbox{\boldmath $E$}(x,T)\cdot\mbox{\boldmath $R$}_e^1)\,dx.
\end{array}
\tag {4.5}
$$
Sine we have
$$\begin{array}{l}
\,\,\,\,\,\,
\displaystyle
\displaystyle
\mbox{\boldmath $\nu$}\times(\mbox{\boldmath $R$}_m^1\times\mbox{\boldmath $\nu$})\cdot
\mbox{\boldmath $\nu$}\times\mbox{\boldmath $R$}_e^1
=-\{\mbox{\boldmath $\nu$}\times(\mbox{\boldmath $R$}_e^1\times\mbox{\boldmath $\nu$})\}\cdot(\mbox{\boldmath $\nu$}\times\mbox{\boldmath $R$}_m^1),
\end{array}
$$
from (3.6) one gets
$$
\displaystyle
\mbox{\boldmath $\nu$}\times(\mbox{\boldmath $R$}_m^1\times\mbox{\boldmath $\nu$})\cdot
\mbox{\boldmath $\nu$}\times\mbox{\boldmath $R$}_e^1
=-\lambda\vert\mbox{\boldmath $\nu$}\times(\mbox{\boldmath $R$}_e^1\times\mbox{\boldmath $\nu$})\vert^2
+\mbox{\boldmath $\nu$}\times(\mbox{\boldmath $R$}_e^1\times\mbox{\boldmath $\nu$})\cdot\mbox{\boldmath $V$}_1.
$$
Thus (4.5) becomes
$$\begin{array}{l}
\displaystyle
\,\,\,\,\,\,
\int_{\Bbb R^3\setminus\overline D}(\tau\mu\vert\mbox{\boldmath $R$}_m^1\vert^2
+\tau\epsilon\vert\mbox{\boldmath $R$}_e^1\vert^2)\,dx
+\int_{\partial D}\lambda\vert\mbox{\boldmath $\nu$}\times(\mbox{\boldmath $R$}_e^1\times\mbox{\boldmath $\nu$})\vert^2\,dS\\
\\
\displaystyle
\,\,\,
+e^{-\tau T}
\int_{\Bbb R^3\setminus\overline D}
(\mu\mbox{\boldmath $H$}(x,T)\cdot\mbox{\boldmath $R$}_m^1+
\epsilon\mbox{\boldmath $E$}(x,T)\cdot\mbox{\boldmath $R$}_e^1)\,dx
\\
\\
\displaystyle
=\int_{\partial D}\mbox{\boldmath $V$}_1\cdot\{\mbox{\boldmath $\nu$}\times(\mbox{\boldmath $R$}_e^1\times\mbox{\boldmath $\nu$})\}
\,dS+F(\tau),
\end{array}
\tag {4.6}
$$
where
$$
\displaystyle
F(\tau)
=-
\int_{\Bbb R^3\setminus\overline D}
\mbox{\boldmath $F$}_1\cdot\mbox{\boldmath $R$}_m^1\,dx
+\int_{\Bbb R^3\setminus\overline D}
\mbox{\boldmath $F$}_2\cdot\mbox{\boldmath $R$}_e^1\,dx.
$$
Rewrite (4.6) further as
$$\begin{array}{l}
\,\,\,\,\,\,
\displaystyle
\int_{\Bbb R^3\setminus\overline D}
\left(
\tau\mu\left\vert
\mbox{\boldmath $R$}_m^1
+\frac{\displaystyle\mbox{\boldmath $F$}_1+e^{-\tau T}\mu\,\mbox{\boldmath $H$}(x,T)}{2\mu\tau}\right\vert^2
+\tau\epsilon\left\vert
\mbox{\boldmath $R$}_e^1
+\frac{\displaystyle -\mbox{\boldmath $F$}_2+e^{-\tau T}\epsilon\mbox{\boldmath $E$}(x,T)}{2\epsilon\tau}\right\vert^2
\right)\,dx\\
\\
\displaystyle
\,\,\,
+\int_{\partial D}\lambda\left\vert\mbox{\boldmath $\nu$}\times(\mbox{\boldmath $R$}_e^1\times\mbox{\boldmath $\nu$})
-\frac{\mbox{\boldmath $V$}_1}{2\lambda}\right\vert^2\,dS\\
\\
\displaystyle
=\int_{\Bbb R^3\setminus\overline D}
\left(
\frac{\displaystyle\vert\mbox{\boldmath $F$}_1+e^{-\tau T}\mu\mbox{\boldmath $H$}(x,T)\vert^2}{4\mu\tau}
+\frac{\displaystyle\vert -\mbox{\boldmath $F$}_2+e^{-\tau T}\epsilon\mbox{\boldmath $E$}(x,T)\vert^2}{4\epsilon\tau}
\right)\,dx
+\int_{\partial D}\frac{\vert\mbox{\boldmath $V$}_1\vert^2}{4\lambda}\,dS.
\end{array}
$$
This immediately yields
$$\begin{array}{l}
\displaystyle
\,\,\,\,\,\,
\Vert\mbox{\boldmath $\nu$}\times(\mbox{\boldmath $R$}_e^1\times\mbox{\boldmath $\nu$})\Vert_{L^2(\partial D)}^2
\\
\\
\displaystyle
\le C(\Vert\mbox{\boldmath $V$}_1\Vert_{L^2(\partial D)}^2
+\tau^{-1}\Vert\mbox{\boldmath $F$}_1\Vert_{L^2(U_{\delta})}^2
+\tau^{-1}\Vert\mbox{\boldmath $F$}_2\Vert_{L^2(U_{\delta})}^2)+O(\tau^{-1}e^{-2\tau T}).
\end{array}
$$
Then, the boundary condition (3.6) yields (4.1).

\noindent
$\Box$

In order to make use of the right-hand side on (4.1), we prepare the following two lemmas.

\proclaim{\noindent Lemma 4.2.}
We have, as $\tau\longrightarrow\infty$
$$
\displaystyle
\lim_{\tau\longrightarrow\infty}
\tau^3 e^{2\tau\sqrt{\mu\epsilon}\,\text{dist}\,(D,B)}
\frac{\displaystyle\Vert\mbox{\boldmath $V$}_1\Vert_{L^2(\partial D)}^2}{\displaystyle\tilde{f}(\tau)^2}
=0.
\tag {4.7}
$$
\endproclaim

{\it\noindent Proof.}
From (2.12) and (3.4) we have the expression
$$\displaystyle
\mbox{\boldmath $V$}_1=-\frac{2\lambda_0\,\lambda}{\lambda+\lambda_0}\mbox{\boldmath $V$}_{em}^0\vert_{\lambda=\lambda_0}.
$$
Thus (2.18) yields
$$\begin{array}{l}
\displaystyle
\,\,\,\,\,\,
\tau^3 e^{2\tau\sqrt{\mu\epsilon}\,\text{dist}\,(D,B)}
\frac{\displaystyle\Vert\mbox{\boldmath $V$}_1\Vert_{L^2(\partial D)}^2}{\displaystyle\tilde{f}(\tau)^2}
\\
\\
\displaystyle
\le C K(\tau)^2\tau^3e^{2\tau\sqrt{\mu\epsilon}\,\text{dist}\,(D,B)}
\int_{\partial D}v^2
\vert\mbox{\boldmath $\nu$}\times({\cal D}(x)\vert_{\lambda=\lambda_0}\mbox{\boldmath $a$}+O(\tau^{-1}))\vert^2\,dS\\
\\
\displaystyle
=C K(\tau)^2\tau^2e^{-2\tau\eta\sqrt{\mu\epsilon}}
\times
\tau e^{2\tau\sqrt{\mu\epsilon}\,d_{\partial D}(p)}
\int_{\partial D}v^2
\vert\mbox{\boldmath $\nu$}\times({\cal D}(x)\vert_{\lambda=\lambda_0}\mbox{\boldmath $a$}+O(\tau^{-1}))\vert^2\,dS.
\end{array}
$$
Note that the term $O(\tau^{-1})$ is uniform with repect to $x\in\partial D$.
Since ${\cal D}(x)\vert_{\lambda=\lambda_0}\mbox{\boldmath $a$}=\mbox{\boldmath $0$}$ for all $x\in\Lambda_{\partial D}(p)$,
it follows from (2.22) that
$$\displaystyle
\lim_{\tau\longrightarrow\infty}\tau
e^{2\tau\sqrt{\mu\epsilon}\,d_{\partial D}(p)}\int_{\partial D}v^2
\vert\mbox{\boldmath $\nu$}\times({\cal D}(x)\vert_{\lambda=\lambda_0}\mbox{\boldmath $a$}+O(\tau^{-1}))\vert^2\,dS
=0.
$$
Then, from (2.23) we obtain the desired conclusion.

\noindent
$\Box$

\proclaim{\noindent Lemma 4.3.}
We have
$$\displaystyle
\Vert\mbox{\boldmath $F$}_1\Vert_{L^2(U_{\delta})}\le C\delta^{-1}\tau^{-1}J_{\infty}(\tau)^{1/2}
\tag {4.8}
$$
and
$$\displaystyle
\Vert\mbox{\boldmath $F$}_2\Vert_{L^2(U_{\delta})}
\le C(\tau^{-1}+\delta+\delta^{-1}\tau^{-1})J_{\infty}^{1/2}(\tau),
\tag {4.9}
$$
where
$$J_{\infty}(\tau)
=\frac{\tau}{\epsilon}\int_{D}\left(\tau\mu\vert\mbox{\boldmath $V$}_m^0\vert^2
+\tau\epsilon\vert\mbox{\boldmath $V$}_e^0\vert^2\right)\,dx.
\tag {4.10}
$$

\endproclaim

{\it\noindent Proof.}
This is an application of a reflection argument developed in \cite{LP}.
First of all, we compute both $\nabla\times(\mbox{\boldmath $V$}_e^0)^*$ and $\nabla\times(\mbox{\boldmath $V$}_m^0)^*$.
From the definition we have
$$\displaystyle
\nabla\times(\mbox{\boldmath $V$}_e^0)^*=-\tau\mu\,(\mbox{\boldmath $V$}_m^0)^*
\tag {4.11}
$$
and hence
$$\displaystyle
\nabla\times(\mbox{\boldmath $V$}_m^0)^*
=-\frac{1}{\tau\mu}\nabla\times\nabla\times(\mbox{\boldmath $V$}_e^0)^*.
\tag {4.12}
$$
Define
$$\displaystyle
\tilde{\phi_{\delta}}(x)
=\frac{\tilde{\lambda}(x)-\lambda_0}{\tilde{\lambda}(x)+\lambda_0}\,\phi_{\delta}(x).
$$
We have
$$\displaystyle
\nabla\times\mbox{\boldmath $R$}_e^0
=\tilde{\phi}_{\delta}\nabla\times(\mbox{\boldmath $V$}_e^0)^*
+\nabla\tilde{\phi}_{\delta}\times(\mbox{\boldmath $V$}_e^0)^*
$$
and from (4.11) one gets the expression
$$\displaystyle
\mbox{\boldmath $F$}_1=\nabla\tilde{\phi}_{\delta}\times(\mbox{\boldmath $V$}_e^0)^*.
$$
This gives
$$\displaystyle
\Vert\mbox{\boldmath $F$}_1\Vert_{L^2(U_{\delta})}
\le C\delta^{-1}\Vert(\mbox{\boldmath $V$}_e^0)^*\Vert_{L^2(U_{\delta})}.
\tag {4.13}
$$
Using the change of variables $y=x^r$, one has
$$
\displaystyle
\Vert(\mbox{\boldmath $V$}_e^0)^*\Vert_{L^2(U_{\delta})}
\le C\Vert\mbox{\boldmath $V$}_e^0\Vert_{L^2(D)}.
\tag {4.14}
$$
We have
$$
\displaystyle
\Vert\mbox{\boldmath $V$}_e^0\Vert_{L^2(D)}
\le \tau^{-1}J_{\infty}(\tau)^{1/2}.
\tag {4.15}
$$
Thus, from (4.13), (4.14) and (4.15), we obtain (4.8).

From (1.9) we know that $\mbox{\boldmath $V$}_e^0$ satisfies
$$\begin{array}{ll}
\displaystyle
\frac{1}{\mu\epsilon}\nabla\times\nabla\times\mbox{\boldmath $V$}_e^0+\tau^2\mbox{\boldmath $V$}_e^0=\mbox{\boldmath $0$}
& \text{in $D$}.
\end{array}
$$
Applying Proposition 3 in \cite{IEM} to this case, we have
$$\begin{array}{c}
\displaystyle
\frac{1}{\mu\epsilon}\nabla\times\nabla\times(\mbox{\boldmath $V$}_e^0)^*+\tau^2(\mbox{\boldmath $V$}_e^0)^*\\
\\
\displaystyle
=\text{terms from $\mbox{\boldmath $V$}_e^0(x^r)$ and $(\mbox{\boldmath $V$}_e^0)'(x^r)$}
+2d_{\partial D}(x)
\times\text{terms from $(\nabla^2\mbox{\boldmath $V$}_e^0)(x^r)$}\\
\\
\displaystyle
\equiv \mbox{\boldmath $Q$}(x)
\end{array}
\tag {4.16}
$$
and all the coefficients in this right-hand side are {\it
independent} of $\tau$ and continuous, in particular, the
coefficients come from the second order terms are $C^1$ in a
tubular neighbourhood of $\partial D$.

From (4.12) we have
$$\begin{array}{ll}
\displaystyle
\nabla\times\mbox{\boldmath $R$}_m^0
&
\displaystyle
=\tilde{\phi}_{\delta}\nabla\times(\mbox{\boldmath $V$}_m^0)^*
+\nabla\tilde{\phi}_{\delta}\times(\mbox{\boldmath $V$}_m^0)^*\\
\\
\displaystyle
&
\displaystyle
=-\frac{\epsilon}{\tau}\tilde{\phi}_{\delta}\,\frac{1}{\mu\epsilon}\nabla\times\nabla\times(\mbox{\boldmath $V$}_e^0)^*
+\nabla\tilde{\phi}_{\delta}\times(\mbox{\boldmath $V$}_m^0)^*.
\end{array}
$$
Thus (4.16) gives
$$\displaystyle
\mbox{\boldmath $F$}_2
=-\frac{\epsilon}{\tau}\tilde{\phi}_{\delta}\mbox{\boldmath $Q$}(x)+\nabla\tilde{\phi}_{\delta}\times(\mbox{\boldmath $V$}_m^0)^*.
$$
This yields
$$\displaystyle
\Vert\mbox{\boldmath $F$}_2\Vert_{L^2(U_{\delta})}
\le
C(
\tau^{-1}\Vert\mbox{\boldmath $Q$}\Vert_{L^2(U_{\delta})}+\delta^{-1}\Vert(\mbox{\boldmath $V$}_m^0)^*\Vert_{L^2(U_{\delta})}).
\tag {4.17}
$$
From the form of $\mbox{\boldmath $Q$}$ and the cahnge of variables, we have
$$\displaystyle
\Vert\mbox{\boldmath $Q$}\Vert_{L^2(U_{\delta})}
\le
C(\Vert\mbox{\boldmath $V$}_e^0\Vert_{L^2(D)}
+\Vert(\mbox{\boldmath $V$}_e^0)'\Vert_{L^2(D)}
+\delta\Vert\nabla^2(\mbox{\boldmath $V$}_e^0)\Vert_{L^2(D)}).
\tag {4.18}
$$
From the definition of $(\mbox{\boldmath $V$}_m^0)^*$ and a change of variables
we have
$$\displaystyle
\Vert(\mbox{\boldmath $V$}_m^0)^*\Vert_{L^2(U_{\delta})}
\le
C\tau^{-1}\Vert(\mbox{\boldmath $V$}_e^0)'\Vert_{L^2(D)}.
\tag {4.19}
$$
Here we claim
$$\displaystyle
\Vert(\mbox{\boldmath $V$}_e^0)'\Vert_{L^2(D)}\le C J_{\infty}(\tau)^{1/2}
\tag {4.20}
$$
and
$$\displaystyle
\Vert\nabla^2(\mbox{\boldmath $V$}_e^0)\Vert_{L^2(D)}\le C \tau J_{\infty}(\tau)^{1/2}
\tag {4.21}
$$
The estimate (4.20) has been established as (27) of Lemma 2.2 in \cite{IEM} since from (2.1) we have another expression
$$\displaystyle
J_{\infty}(\tau)=\frac{1}{\mu\epsilon}\int_D\vert\nabla\times\mbox{\boldmath $V$}_e^0\vert^2\,dx
+\tau^2\int_D\vert\mbox{\boldmath $V$}_e^0\vert^2\,dx.
$$

The estimate (4.21) is proved using the explicit from (2.16) as follows.
We have
$$\begin{array}{l}
\displaystyle
\,\,\,\,\,\,
(\mbox{\boldmath $V$}_e^0)'(x)
\\
\\
\displaystyle
=K(\tau)\tilde{f}(\tau)\left\{v(x)(\mbox{\boldmath $M$}(x;p)\mbox{\boldmath $a$})'
+(\mbox{\boldmath $M$}(x;p)\mbox{\boldmath $a$})\otimes\nabla v(x)\right\}\\
\\
\displaystyle
=K(\tau)\tilde{f}(\tau)v(x)\{(\mbox{\boldmath $M$}(x;p)\mbox{\boldmath $a$})'-\tilde{\tau}
\left(1+\frac{1}{\tilde{\tau}\vert x-p\vert}\right)
(\mbox{\boldmath $M$}(x;p)\mbox{\boldmath $a$})\otimes\mbox{\boldmath $\omega$}_x\}
\end{array}
$$
where $\tilde{\tau}=\tau\,\sqrt{\mu\epsilon}$ and
$$
\begin{array}{l}
\displaystyle
\,\,\,\,\,\,
(\mbox{\boldmath $M$}(x;p)\mbox{\boldmath $a$})'\\
\\
\displaystyle
=\frac{1}{\tilde{\tau}}\left(\frac{1}{\vert x-p\vert^2}+\frac{2}{\tilde{\tau}\vert x-p\vert^3}\right)
\left(3\mbox{\boldmath $\omega$}_x\otimes\mbox{\boldmath $\omega$}_x\,\mbox{\boldmath $\omega$}_x\cdot\mbox{\boldmath $a$}
-\mbox{\boldmath $a$}\otimes\mbox{\boldmath $\omega$}_x\right)\\
\\
\displaystyle
\,\,\,
-
\left\{\frac{1}{\vert x-p\vert}+\frac{3}{\tilde{\tau}}
\left(\frac{1}{\vert x-p\vert^2}+\frac{1}{\tilde{\tau}\vert x-p\vert^3}\right)
\right\}\left(
\mbox{\boldmath $\omega$}_x\cdot\mbox{\boldmath $a$}
I_3
-2\mbox{\boldmath $\omega$}_x\cdot\mbox{\boldmath $a$}\,\mbox{\boldmath $\omega$}_x
\otimes
\mbox{\boldmath $\omega$}_x
+\mbox{\boldmath $\omega$}_x
\otimes
\mbox{\boldmath $a$}
\right).
\end{array}
$$
Then, from (2.16) we see that
$$\begin{array}{l}
\displaystyle
(\mbox{\boldmath $V$}_e^0)'(x)
=-\tilde{\tau}\left(1+\frac{1}{\tilde{\tau}\vert x-p\vert}\right)\mbox{\boldmath $V$}_e^0(x)\otimes\mbox{\boldmath $\omega$}_x
+K(\tau)\tilde{f}(\tau)v(x)(\mbox{\boldmath $M$}(x;p)\mbox{\boldmath $a$})'
\end{array}
$$
and hence, for $j=1,2,3$
$$\begin{array}{l}
\displaystyle
\,\,\,\,\,\,
\frac{\partial}{\partial x_j}(\mbox{\boldmath $V$}_e^0)'(x)
\\
\\
\displaystyle
=\frac{(x_j-p_j)}{\vert x-p\vert^3}\mbox{\boldmath $V$}_e^0(x)\otimes\mbox{\boldmath $\omega$}_x
-\tilde{\tau}\left(1+\frac{1}{\tilde{\tau}\vert x-p\vert}\right)\{\frac{\partial}{\partial x_j}\mbox{\boldmath $V$}_e^0(x)\}\otimes\mbox{\boldmath $\omega$}_x
\\
\\
\displaystyle
-\tilde{\tau}\left(1+\frac{1}{\tilde{\tau}\vert x-p\vert}\right)\mbox{\boldmath $V$}_e^0(x)\}\otimes\frac{\partial}{\partial x_j}\mbox{\boldmath $\omega$}_x
\\
\\
\displaystyle
+K(\tau)\tilde{f}(\tau)\frac{\partial}{\partial x_j}v(x)(\mbox{\boldmath $M$}(x;p)\mbox{\boldmath $a$})'
+K(\tau)\tilde{f}(\tau)v(x)\frac{\partial}{\partial x_j}(\mbox{\boldmath $M$}(x;p)\mbox{\boldmath $a$})'.
\end{array}
$$
Then, it is easy to see that, there exists a positive constant $C$ independent of $\tau$ such that, for all $x\in D$ and $\tau>0$, we have
$$\displaystyle
\vert\nabla^2(\mbox{\boldmath $V$}_e^0)(x)\vert
\le C(\tau+1)(\vert\mbox{\boldmath $V$}_e^0(x)\vert+\vert(\mbox{\boldmath $V$}_e^0)'(x)\vert+K(\tau)\vert\tilde{f}(\tau)\vert\vert v(x)\vert).
$$
We know from (24) in \cite{IEM} that, for all $x\in\Bbb R^3\setminus B$
$$\displaystyle
\vert\mbox{\boldmath $V$}_e^0(x)\vert^2
\ge \tau^{-2}\frac{\displaystyle K(\tau)^2\tilde{f}(\tau)^2 v(x)^2}{\displaystyle \mu\epsilon\vert x-p\vert^2}.
$$
This yields
$$\displaystyle
\Vert\nabla^2(\mbox{\boldmath $V$}_e^0)\Vert_{L^2(D)}^2
\le C'\tau^2(\Vert\mbox{\boldmath $V$}_e^0\Vert_{L^2(D)}^2+
\Vert(\mbox{\boldmath $V$}_e^0)'\Vert_{L^2(D)}^2).
$$
Then, applying (4.15) and (4.20) to this right-hand side we obtain (4.21).

Now applying (4.15), estimates (4.20) and (4.21) to (4.18), one gets
$$\displaystyle
\Vert\mbox{\boldmath $Q$}\Vert_{L^2(U_{\delta})}
\le C(1+\delta\tau)J_{\infty}(\tau)^{1/2}.
$$
Then, this together with (4.17), (4.19) and (4.20) yields (4.9).

\noindent
$\Box$

A combination of (2.2) in $D$ and (4.10), we have
$$\begin{array}{ll}
\displaystyle
J_{\infty}(\tau)
&
\displaystyle
=-\frac{\tau}{\epsilon}\int_{\partial D}\mbox{\boldmath $\nu$}\cdot(\mbox{\boldmath $V$}_e^0\times\mbox{\boldmath $V$}_m^0)\,dS\\
\\
\displaystyle
&
\displaystyle
=\frac{\tau}{\epsilon}
\int_{\partial D}
\mbox{\boldmath $\nu$}\times(\mbox{\boldmath $V$}_e^0\times\mbox{\boldmath $\nu$})\cdot\mbox{\boldmath $\nu$}\times\mbox{\boldmath $V$}_m^0\,dS.
\end{array}
$$
Since this last expreesion means that $J_{\infty}(\tau)=J(\tau)\vert_{\lambda=\infty}$, similarly to (2.13), we have
$$\begin{array}{l}
\displaystyle
\,\,\,\,\,\,
\lim_{\tau\longrightarrow\infty}
\tau^2 e^{2\tau\sqrt{\mu\epsilon}\,\text{dist}\,(D,B)}\frac{J_{\infty}(\tau)}{\displaystyle\tilde{f}(\tau)^2}\\
\\
\displaystyle
=\frac{\pi}{4}\left(\frac{\eta}{d_{\partial D}(p)}\right)^2\frac{\lambda_0^2}{\epsilon^4}
\sum_{q\in\Lambda_{\partial D}(p)}k_q(p)
\vert\mbox{\boldmath $\nu$}_q\times(\mbox{\boldmath $a$}\times\mbox{\boldmath $\nu$}_q)\vert^2.
\end{array}
$$
This gives, as $\tau\longrightarrow\infty$
$$\displaystyle
e^{2\tau\sqrt{\mu\epsilon}\,\text{dist}\,(D,B)}\frac{J_{\infty}(\tau)}{\displaystyle\tilde{f}(\tau)^2}=O(\tau^{-2}).
$$
From this together with (4.8) and (4.1), we have
$$\begin{array}{l}
\displaystyle
\,\,\,\,\,\,
\tau^{-1}\Vert\mbox{\boldmath $F$}_1\Vert_{L^2(U_{\delta})}^2
+\tau^{-1}\Vert\mbox{\boldmath $F$}_2\Vert_{L^2(U_{\delta})}^2\\
\\
\displaystyle
\le
C(\delta^{-2}\tau^{-5}+\tau^{-5}+\delta^2\tau^{-3})e^{-2\tau\sqrt{\mu\epsilon}\,\text{dist}\,(D,B)}\tilde{f}(\tau)^2
\end{array}
$$
and hence
$$\begin{array}{l}
\displaystyle
\,\,\,\,\,\,
\tau^3(\tau^{-1}\Vert\mbox{\boldmath $F$}_1\Vert_{L^2(U_{\delta})}^2
+\tau^{-1}\Vert\mbox{\boldmath $F$}_2\Vert_{L^2(U_{\delta})}^2)\\
\\
\displaystyle
\le
C(\delta^{-2}\tau^{-2}+\tau^{-2}+\delta^2)e^{-2\tau\sqrt{\mu\epsilon}\,\text{dist}\,(D,B)}\tilde{f}(\tau)^2.
\end{array}
\tag {4.22}
$$
Now choosing $\delta=\tau^{-\theta}$ with $\theta>0$, we have
$\delta^{-2}\tau^{-2}+\tau^{-2}+\delta^2=\tau^{-2(1-\theta)}+\tau^{-2}+\tau^{-2\theta}$.
Thus choosing $\theta=1/2$, we have $1-\theta=\theta$ and (4.22) becomes
$$
\displaystyle
\tau^3(\tau^{-1}\Vert\mbox{\boldmath $F$}_1\Vert_{L^2(U_{\delta})}^2
+\tau^{-1}\Vert\mbox{\boldmath $F$}_2\Vert_{L^2(U_{\delta})}^2)
\le
C\tau^{-1}e^{-2\tau\sqrt{\mu\epsilon}\,\text{dist}\,(D,B)}\tilde{f}(\tau)^2.
$$
Now applying this and (4.7) to the right-hand side on (4.1) together with (1.8) and (2.9),
we obtain (3.13).  This completes the proof of Lemma 3.2.

$$\quad$$

\centerline{{\bf Acknowledgments}}

The author was partially supported by Grant-in-Aid for
Scientific Research (C)(No 17K05331) of Japan Society for
the Promotion of Science.

$$\quad$$

\vskip1cm
\noindent
e-mail address

ikehata@hiroshima-u.ac.jp

\end{document}